\documentclass[english, hidelinks, 11pt]{article}
\usepackage{geometry}
\usepackage{setspace}
\usepackage{bm}
\usepackage{graphicx}
\graphicspath{ {.} }
\usepackage{xr}

\usepackage{amssymb, amsmath, amsthm, amsfonts, mathtools, bbm, mathrsfs}
\usepackage{float}
\usepackage{color}
\usepackage[dvipsnames]{xcolor}
\usepackage{url}
\usepackage{enumitem}
\usepackage{csquotes}
\usepackage{multirow}
\usepackage{siunitx}
\usepackage{placeins}
\usepackage{authblk}
\usepackage{appendix}
\usepackage{lscape}
\usepackage{booktabs}
\usepackage{longtable}
\usepackage{stackrel}

\usepackage{tabularx}

\allowdisplaybreaks[1]

\usepackage{tikz}
\usetikzlibrary{arrows.meta, decorations.pathreplacing}

\usepackage{algorithm}
\usepackage{algpseudocode}

\usepackage{hyperref}
\hypersetup{
    colorlinks=true,
    linkcolor=blue,
    citecolor=blue,
    filecolor=magenta,
    urlcolor=blue
}
\usepackage{cleveref}

\usepackage[round]{natbib}

\newtheorem{theorem}{Theorem}
\newtheorem{lemma}{Lemma}
\newtheorem{proposition}{Proposition}
\newtheorem{corollary}{Corollary}

\newtheorem{remark}{Remark}
\newtheorem{assumption}{Assumption}

\numberwithin{equation}{section}

\usepackage{soul}
\usepackage{ulem}
\DeclareMathOperator*{\argmax}{arg\,max}

\DeclarePairedDelimiter{\floor}{\lfloor}{\rfloor}

\def\tr{\textnormal{\text{tr}}}

\def\wh#1{\widehat{#1}}

\def\E{\mathbb{E}}

\def\P{\mathbb{P}}
\def\R{\mathbb{R}}

\def\S{\mathbb{S}}

\def\cB{\mathcal{B}}

\def\cF{\mathcal{F}}

\def\cI{\mathcal{I}}

\def\cM{\mathcal{M}}
\def\cN{\mathcal{N}}
\def\cO{\mathcal{O}}

\def\cT{\mathcal{T}}

\def\zero{\mathbf{0}}

\def\ba{\mathbf{a}}
\def\bb{\mathbf{b}}

\def\bs{\mathbf{s}}
\def\be{\mathbf{e}}
\def\bn{\mathbf{n}}
\def\bff{\mathbf{f}}
\def\bu{\mathbf{u}}
\def\bh{\mathbf{h}}
\def\bx{\mathbf{x}}
\def\by{\mathbf{y}}
\def\bv{\mathbf{v}}
\def\bw{\mathbf{w}}
\def\bg{\mathbf{g}}
\def\bq{\mathbf{q}}
\def\bf{\mathbf{f}}

\def\bz{\mathbf{z}}
\def\bA{\mathbf{A}}

\def\bD{\mathbf{D}}

\def\bG{\mathbf{G}}
\def\bH{\mathbf{H}}
\def\bI{\mathbf{I}}

\def\bM{\mathbf{M}}

\def\bP{\mathbf{P}}
\def\bQ{\mathbf{Q}}
\def\bR{\mathbf{R}}
\def\bS{\mathbf{S}}
\def\bL{\mathbf{L}}
\def\bU{\mathbf{U}}
\def\bV{\mathbf{V}}

\def\bX{\mathbf{X}}
\def\bY{\mathbf{Y}}

\def\bmu{\boldsymbol{\mu}} \def\bLambda{\boldsymbol{\Lambda}}

 \def\bSigma{\boldsymbol{\Sigma}}
 
 \def\bxi{\boldsymbol{\xi}}

\def\bPi{\boldsymbol{\Pi}} \def\bepsilon{\boldsymbol{\epsilon}} 
\def\bchi{\boldsymbol{\chi}}

\def\spa{\mathrm{span}}

\def\diag{\mathrm{diag}}

\definecolor{hancolor}{RGB}{170,74,20}

\usepackage[makeroom]{cancel}

\usepackage{parskip}

\begin{document}
\def\spacingset#1{\renewcommand{\baselinestretch}%
{#1}\small\normalsize} \spacingset{1}

	\begin{titlepage}

        \title{Divide-and-Conquer is Tail- and Anomaly-Robust: Factor Model Estimation Without Moment Constraints}
                \author{Zetai Cen\thanks{Corresponding author. School of Mathematics, University of Bristol.
                Email: \url{zetai.cen@bristol.ac.uk}.
                Supported by Engineering and Physical Sciences Research Council (EP/Z531327/1).}}
                \author{Han Yan\thanks{Department of Statistics, London School of Economics.
                Email: \url{H.Yan22@lse.ac.uk}.
                Supported by Engineering and Physical Sciences Research Council (EP/Z531327/1).}}
                \author{Haeran Cho\thanks{School of Mathematics, University of Bristol.
                Email: \url{haeran.cho@bristol.ac.uk}.
                Supported by Engineering and Physical Sciences Research Council (EP/Z531327/1).}}

        \affil{School of Mathematics, University of Bristol}
        \affil{Department of Statistics, London School of Economics}

		\date{}

		\maketitle

\begin{abstract}
We study robust estimation of elliptical factor models with abnormal observations and without requiring the existence of any positive moments.
In particular, a divide-and-conquer (DC) estimator averages orthogonal projectors onto the full column spaces of small data blocks,
and we rigorously establish its probability bound that allows weak and heterogeneous factors.
The derived convergence rate is comparable to that of classical principal component estimators, up to the proportion of anomalies.
Building on our DC estimator, we propose a split-sample tail- and anomaly-robust (STAR) estimator for the factor number using sample medians of successive projected coordinate ratios,
which is shown consistent under mild conditions.
Subsequent factor and common component estimation is also discussed, with the rates of convergence spelled out.
Remarkably, given an upper bound on the true factor number, our suite of procedures are free from tuning,
require no tail-dependent truncation threshold nor solving an optimization problem,
admit parallel block computations,
and hence has the same computational complexity as the classical principal component methods.
Extensive experiments corroborate our theoretical findings, where DC and STAR are comparable and often outperform state-of-the-art competitors.
We also illustrate our methods on a portfolio return dataset.
\end{abstract}
		
		\bigskip
		\bigskip


\noindent

	\end{titlepage}

	\setcounter{page}{2}

\tableofcontents

\renewcommand{\baselinestretch}{1.8}

\newpage
\section{Introduction}\label{sec: introduction}

High-dimensional factor models describe common variation in many observed variables through a small number of latent factors.
They provide a framework for dimension reduction, forecasting, and the analysis of economic and financial panels.
Principal component methods and their large-dimensional theory have played a central role in these developments \citep{StockWatson2002, BaiNg2002, Bai2003, barigozzi2026principal}, where the primary focus lies in estimating the vector space spanned by the matrix of loadings, a.k.a.\ \textit{loading space}.
Another critical part of factor analysis is to determine the number of factors, and several approaches have been proposed which are based on information criteria or eigenvalue ratios \citep{BaiNg2002, ahn2013eigenvalue, LamYao2012}.
While it is commonly assumed that the factors are pervasive across the cross-sections, there also exist studies for understanding how the strength of factors affects the rates of estimation \citep{Lametal2011, BaiNg2023}.

Two distinct stylized features of modern data complicate the use of conventional principal component methods.
First, the data may have heavy tails as frequently observed in economics and finance \citep{cont2001}, which invalidates the commonly made moment conditions underpinning covariance-based approaches.
Second, some observations may be anomalous and arbitrarily deviate from the underlying model, due to the events such as financial crisis.
When heavy tails and anomalies co-exist, extreme observations may be the realizations from the heavy-tailed distribution or anomalies and, regardless of their origin, negatively influence the accuracy of the existing principal component-based estimators \citep{Hubertetal2005, HeKongYuZhang2022}.
In this paper, we develop methods for estimating the loading space and the number of factors under elliptical factor models while permitting the presence of heavy tails, anomalies, and weak factors. 

Elliptical distributions are popularly considered for high-dimensional data analysis under heavy tails.
For instance, \cite{HanLiu2018} develop elliptical component analysis using multivariate Kendall's tau, and \cite{FanLiuWang2018} study robust covariance estimation through elliptical factor models.
For factor modeling, methods based on spatial Kendall's tau have been developed for estimating the loading space and the latent factor process \citep{HeKongYuZhang2022} and the factor number \citep{Yuetal2019}.
Other methods address heavy tails under different distributional assumptions, e.g.,\ by adopting the Spearman correlation for factor number estimation under a scale-mixture structure \citep{Qiuetal2025}, utilizing Huber loss \citep{Barigozzietal2023_RTFA} or data truncation \citep{Barigozzietal2026_robust}.
These approaches offer different tradeoffs between distributional structure and moment requirements.

There also exists a closely related but distinct line of research which focuses on developing principal component methods that are robust against outliers.
Such methods combine the projection pursuit with robust scatter estimation \citep{Hubertetal2005}, and can accommodate sparse loadings \citep{Hubertetal2016}, or cellwise and rowwise outliers \citep{Hubertetal2019}. 
It is worth pointing out that the procedures designed to control heavy-tailed fluctuations and those designed to be robust against the anomalies, address related but different objectives; see \cite{raymaekers2026challenges} for a comprehensive overview and further discussions.

In this paper, we propose the \textbf{D}ivide-and-\textbf{C}onquer (DC) estimator for the loading space under elliptical factor models.
Simply put, DC partitions the data into small blocks, computes the orthogonal projector onto the full column space of each block, and extracts the leading eigenspace of their average.
When the factor number $r$ is known, we take the block size to be $r$; if it is unknown, we set the block size to be a strict upper bound on $r$. 
This sets the DC estimator apart from the existing tail-robust methods that often rely on selecting the tuning parameters that depend on the unknown tail index.
The aggregation of local spectral information has an established role in distributed estimation \citep{Fanetal2019, GaoTsay2023, HeLiuWang2025} but in contrast to these earlier works, the DC estimator aggregates the blockwise orthogonal projectors of full column rank of the block size,
even when the block size exceeds the number of factors.
Our investigation reveals that this proposed approach allows for establishing a signal-to-noise eigengap under possibly weak and heterogeneous factor strengths, without making any moment constraint.

For factor number estimation, we propose the \textbf{S}plit-sample \textbf{T}ail- and \textbf{A}nomaly-\textbf{R}obust (STAR) estimator. 
splitting the data into two, the STAR projects the second half of the data onto the DC estimator from the first half, takes the ratio of each pair of successive squared scores, and locates the maximizer of the sample median of these ratios over time.
By its construction, the STAR estimator utilizes the ellipticity of the data and benefits from the eigengap preserved in the DC loading space estimator.


Theoretically, without imposing any finite-moment requirements while permitting the presence of anomalies,
we establish a non-asymptotic probabilistic bound for the estimation error of the DC estimator.
If factors are strong and the anomalies are absent,
by setting the block size $b$ as the factor number ($b = r$), the convergence rate of the DC estimator
matches the rate attained by the classical principal component methods under the existence of the fourth moment in the data \citep{BaiNg2002}. 
We also show the consistency of the STAR estimator with which the additional application of the DC estimator allows us to estimate the factors and common components, also without any moment constraints.

To summarize, Tables~\ref{tab: loading_comparison} and \ref{tab: rank_comparison} provide an overview of the methods most closely related to our study in comparison with DC and STAR estimators, respectively.
For \cite{Barigozzietal2026_robust}, the displayed rate concerns the vector factor model under their general dependence conditions. Under their bounded deterministic factors and the remaining conditions of their loading estimation result, their rate would reduce to $\cO_P\big( 1/\sqrt{T} + 1/p \big)$.

\begin{table}[th!]
\centering
\small
\setlength{\tabcolsep}{3pt}
\renewcommand{\arraystretch}{1.25}
\caption{
Comparison of loading space and principal subspace estimators.
The rates are derived for the distance between true loading matrix $\bA$ and estimator $\wh\bA$ measured by $\big\| \bA (\bA^\top \bA)^{-1} \bA^\top - \wh\bA (\wh\bA^\top \wh\bA)^{-1} \wh\bA^\top \big\|_F$, 
under strong factors and in the absence of anomalies; the rate reported for DC is a special case of Theorem~\ref{thm: DC_ano_clean}. Here, $T$ refers to the sample size and $p$ the dimension of the data.
}
\label{tab: loading_comparison}

\begin{tabularx}{\textwidth}{
    @{}
    >{\raggedright\arraybackslash}p{0.22\textwidth}
    >{\raggedright\arraybackslash}p{0.20\textwidth}
    >{\raggedright\arraybackslash}p{0.24\textwidth}
    >{\raggedright\arraybackslash}X
    @{}
}
\toprule
Method
& Moment requirement
& Model
& Convergence rate \\
\midrule

\cite{Bai2003}
& $\geq 4$
& General factor model
& $\displaystyle
   \cO_P\bigg(
   \sqrt{\frac{1}{T}}
   + \frac{1}{p}
   \bigg)$
\\

\cite{Barigozzietal2026_robust}
& $2 + 2\epsilon$, $\epsilon>0$
& General factor model
& $\displaystyle
   \cO_P\bigg(
   \bigg\{\frac{\log(Tp)}{T} \bigg\}^{\epsilon/(1+\epsilon)}
   + \frac{1}{p}
   \bigg)$
\\

\cite{HanLiu2018}
& None
& Elliptical factor model
& $\displaystyle
   \cO_P\bigg(
   \sqrt{\frac{\log(p)}{T}}
   \bigg)$
\\

\cite{HeKongYuZhang2022}
& None
& Elliptical factor model
& $\displaystyle
   \cO_P\bigg(
   \sqrt{\frac{1}{T}}
   + \frac{1}{p}
   \bigg)$
\\


\cite{Hubertetal2005}
& Unspecified
& Robust PCA
& N/A
\\

\cite{Hubertetal2016}
& Unspecified
& Robust PCA
& N/A
\\

\cite{Hubertetal2019}
& Unspecified
& Robust PCA
& N/A 
\\

DC
& None
& Elliptical factor model
& $\displaystyle
   \cO_P\bigg(
   \sqrt{\frac{1}{T}}
   + \frac{1}{p}
   \bigg)$
\\

\bottomrule
\end{tabularx}

\vspace{0.6em}

\caption{
Comparison of factor number estimators.
}
\label{tab: rank_comparison}

\begin{tabularx}{\textwidth}{
    @{}
    >{\raggedright\arraybackslash}p{0.35\textwidth}
    >{\raggedright\arraybackslash}p{0.23\textwidth}
    >{\raggedright\arraybackslash}X
    @{}
}
\toprule
Method
& Moment requirement
& Model \\
\midrule

\cite{BaiNg2002}
& $\geq 4$ 
& General factor model
\\

\cite{Barigozzietal2026_robust}
& $2 + 2\epsilon$, $\epsilon>0$
& General factor model
\\

\cite{FanGuoZheng2022}
& $6 + \epsilon$, $\epsilon>0$
& Independent latent coordinates
\\

\cite{Hubertetal2005}
& Unspecified
& Robust PCA
\\

\cite{Qiuetal2025}
& None
& Scale-mixture factor model
\\


\cite{Yuetal2019}
& None
& Elliptical factor model
\\

STAR
& None
& Elliptical factor model
\\

\bottomrule
\end{tabularx}

\end{table}


\textbf{Notations.}
Unless specified otherwise, we denote scalars, vectors and matrices by plain letters, lower-case bold letters and capital bold letters, i.e., $x$ (or $X$), $\bx$ and $\bX$, respectively. 
We use $x_i$ and $X_{ij}$ to denote the $i$th element of a vector $\bx$ and the $(i,j)$th element of a matrix $\bX$, respectively.
We use $a\asymp b$ (\text{resp.} $a\asymp_P b$) to denote $a=\cO(b)$ and $b=\cO(a)$ (\text{resp.} $a=\cO_P(b)$ and $b=\cO_P(a)$). 
Hereafter, given a positive integer $m$, define $[m]:=\{1,2,\dots,m\}$.
The $i$th largest eigenvalue of a real symmetric matrix $\bX$ is denoted by $\lambda_i(\bX)$,
and $\mathrm{SVD}_i(\bX)$ represents the leading $i$ eigenvectors of $\bX$, by which
we always refer to the orthonormal eigenvectors of the matrix corresponding to the $i$ largest eigenvalues.
For any symmetric matrices $\bX_1$ and $\bX_2$, the notation $\bX_1 \preceq \bX_2$ means $\bX_2-\bX_1$ is positive semi-definite.

\section{Methodology}\label{sec: method}

\subsection{Model and problem formulation}\label{subsec: model}

Consider a sequence of $p$-dimensional random vectors $\{\bx_t^0\}_{t\in[T]}$
admitting the factor model representation \citep[e.g.][]{BaiNg2002, Bai2003} such that
\begin{equation}\label{eqn: vfm}
\bx_t^0 = \bA \bff_t + \bepsilon_t =: \bchi_t + \bepsilon_t ,
\end{equation}
where $\bA \in \R^{p \times r}$ is the factor loading matrix,
$\bff_t \in \R^r$ is the latent core factor,
and $\bepsilon_t \in \R^p$ is the idiosyncratic noise.
The factor-driven $\bchi_t$ is referred to as the common component.
The number of factors $r$ is fixed but unknown.
Then, we observe $\{\bx_t\}_{t\in[T]}$, where
\begin{align}
\label{eq:obs}
\bx_t =
\begin{cases}
\bx_t^0 , & t \notin \cT , \\
\bz_t , & t \in \cT ,
\end{cases}
\end{align}
for an unknown set $\cT \subset [T]$ and arbitrary vectors $\bz_t \in \R^p$ representing the anomalies. 
The index set of anomalies may be random or deterministic, see Assumption~\ref{ass: contamination} below for the conditions on $\cT$.
In this setup, it is of interest to estimate the \textit{loading space} spanned by the columns of $\bA$ and the factor number $r$ under~\eqref{eqn: vfm}, from the observations $\{\bx_t\}_{t\in[T]}$ which are potentially contaminated by anomalies. 
Further, even in the absence of anomalies, we permit $\{\bx^0_t\}_{t\in[T]}$ to be heavy-tailed without any moment condition as specified below.

Following \cite{FanLiuWang2018} and \cite{HeKongYuZhang2022}, we formulate the model for $\bx_t^0$ through a joint elliptical distribution.
We say a random vector $\by \in \R^{h}$ follows an elliptical distribution, denoted by $\by \sim \mathrm{ED}_h(\bmu, \bSigma_0, \zeta)$, if
\begin{equation}
\label{eqn: elliptical}
\by \overset{d}{=} \bmu + \zeta \bA_0 \bu ,
\end{equation}
where $\bmu \in \R^{h}$ is deterministic, $\bu \in\R^{q}$ is a random vector uniformly distributed on the unit sphere $\S^{q-1}$ with $q\leq h$, 
the radius $\zeta \geq 0$ is a scalar random variable independent of $\bu$, $\bA_0 \in \R^{h \times q}$ is a deterministic matrix with full column rank, and the matrix $\bSigma_0 = \bA_0 \bA_0^\top$ is referred to as the scatter matrix of $\by$.
Such framework contains multivariate Gaussian distribution as a special case and, more importantly, includes a large class of heavy-tailed distributions such as the multivariate elliptical $t$-distribution.

Note that the scatter matrix $\bSigma_0$ is only defined up to a positive scalar factor.
Namely, the distribution in~\eqref{eqn: elliptical} would remain unchanged if $\bA_0$ and $\zeta$ are replaced by $c \bA_0$ and $c^{-1} \zeta$, respectively.
In view of this, in formulating the elliptical distribution of $\bx^0_t$, we adopt a representation where the block of the joint scatter matrix corresponding to the factor is $\bI_r$.

\begin{assumption}\label{ass: elliptical}
For each $t\in[T]$, we assume that $\bff_t$ and $\bepsilon_t$ in~\eqref{eqn: vfm} 
satisfy
\begin{equation}\label{eqn: elliptical_obs}
\begin{pmatrix} \bff_t \\ \bepsilon_t \end{pmatrix}
=
\zeta_t
\begin{pmatrix} \bI_r & \zero \\ \zero & \bA_\epsilon \end{pmatrix}
\frac{\bu_t}{\|\bu_t\|} ,
\end{equation}
where $\zeta_t$ and $\bu_t$ are independent copies of a scalar
random variable $\zeta$ and a random vector $\bu\sim \cN(0,\bI_{r+p})$, respectively.
Also, $\zeta$ is independent of $\bu$,
and is finite and strictly positive almost surely.
\end{assumption}

Assumption~\ref{ass: elliptical} implies that
$(\bff_t^\top, \bepsilon_t^\top)^\top \sim \mathrm{ED}_{r+p}(\zero, \diag(\bI_r, \bSigma_\epsilon), \zeta)$, where $\bSigma_\epsilon := \bA_\epsilon \bA_\epsilon^\top \in \R^{p \times p}$.
Since linear transformations of elliptically distributed random vectors remain elliptical, 
the model in~\eqref{eqn: vfm} further implies that $\bx_t^0$ is elliptically distributed itself.
Specifically, let us 
define $\bM := [\bA, \bA_\epsilon] \in \R^{p \times (r+p)}$ and $\bv_t := \bu_t / \|\bu_t\|$ for $t\in[T]$.
Then $\bx_t^0$ is elliptically distributed and admits the representation
\begin{equation}
\label{eqn: xt_ellip_form}
\bx_t^0 = \zeta_t \bM \bv_t .
\end{equation}
with its scatter matrix 
\begin{equation}\label{eqn: cov}
\bSigma_x := 
\bM \bM^\top
= \bA \bA^\top + \bSigma_\epsilon .
\end{equation}

Without loss of generality, we suppose that $\zeta$ is nonnegative as the sign can be absorbed into the direction vector $\bu_t$.
Such convention is standard in the stochastic representation of elliptical distributions \citep{FanLiuWang2018, HeKongYuZhang2022}.
Assumption~\ref{ass: elliptical} particularly requires $\P(\zeta >0) = 1$, which can be seen as a nondegeneracy condition.
Also, we do not place any restriction on the growth of $\zeta_t$ with $p$, unlike in \cite{HeKongYuZhang2022} where it is required that $\zeta/\sqrt{p} = \cO_P(1)$, see their Assumption~A.
Further, we do not impose any moment condition on $\bx_t^0$, viz.\ $\bx_t^0$ is not required to have finite first or second moments, and the scatter matrix $\bSigma_x$ is a well-defined quantity that characterizes the cross-sectional dependence of the data even when the covariance matrix of $\bx_t^0$ 
does not exist.


\begin{assumption}
\label{ass: loading}
The number of factors $r$ is fixed.
For some fixed $\delta > 0$ and an absolute constant $c_\lambda > 0$,
the loading matrix satisfies $\bA^\top \bA \succeq c_\lambda p^{\delta} \bI_r$
for all sufficiently large $p$.
\end{assumption}

\begin{assumption}\label{ass: noise}
The matrix $\bSigma_\epsilon$ satisfies $c_\epsilon \bI_p \preceq \bSigma_\epsilon \preceq C_\epsilon \bI_p$ for some absolute constants $c_\epsilon, C_\epsilon >0$.
\end{assumption}



Assumption~\ref{ass: loading} provides a very general setup, in that it only imposes a lower bound on the weakest loading direction, and permits the factor strengths captured by the eigenvalues of $\bA^\top\bA$ to be highly heterogeneous.
The case when $0 < \delta < 1$ allows for weak factors, while $\delta = 1$ indicates that the factors are pervasive across the cross-sections \citep[cf.][Section~2.3]{Lametal2011}, and we allow for the case where $\delta > 1$.
Assumption~\ref{ass: noise} bounds the noise scatter eigenvalues away from zero and infinity.
From~\eqref{eqn: cov}, Assumptions~\ref{ass: loading} and \ref{ass: noise}, we have $\lambda_r(\bSigma_x) \geq c_\lambda p^{\delta} + c_\epsilon$ for all sufficiently large $p$, while $c_\epsilon \leq \lambda_j(\bSigma_x) \leq C_\epsilon$ for $r < j \leq p$,
so that the leading $r$ eigenvalues of $\bSigma_x$ are separated from the remaining ones.

We allow a subset of the observations to be replaced by arbitrary anomalous vectors, as formalized below.

\begin{assumption}\label{ass: contamination}
There exists a constant
$\overline{\cT} \in [0, T)$
such that the possibly random anomaly index set $\cT$ is independent of $\{(\zeta_t, \bu_t)\}_{t\in [T]}$ in Assumption~\ref{ass: elliptical}, and satisfies
\begin{equation}
\label{eqn: contam_index}
\cT \subset [T] \text{ \ and \ }
|\cT| \leq \overline{\cT}
\quad \text{almost surely.}
\end{equation}
\end{assumption}

It is worth pointing out that the bound $\overline{\cT}$ in Assumption~\ref{ass: contamination} may depend on $T$ and $p$, and the anomalous vectors may depend on the latent $\{\bx^0_t\}_{t \in [T]}$ and need not have bounded magnitudes nor belong to any specified subspace.
Throughout, the set $\cT$ and the anomalous observations are unknown.


\subsection{Divide-and-Conquer estimation of the loading space}

We describe the proposed DC loading space estimation method, which proceeds by combining blockwise estimators of the loading space from the data blocks of size $b$.
We first present the method when the factor number $r$ is known, in which case we set the block size as $b = r$.
More generally, the following construction is defined for any integer $r \leq b \leq \min(T,p)$.

Set $m := \floor{T/b}$ and partition the first $mb$ observations into consecutive blocks $\cM_i := \{(i-1)b + 1, \dots, ib\}$, $i\in[m]$.
For each block~$i\in[m]$, compute the following matrix and thus its $b$ leading eigenvectors as
\begin{equation}
\label{eqn: block_cov_vfm}
\wh\bS_i := \frac{1}{b} \sum_{t \in \cM_i} \bx_t \bx_t^\top 
\in \R^{p \times p},
\qquad
\wh\bU_i := \mathrm{SVD}_{b}(\wh\bS_i) 
\in \R^{p \times b} ,
\end{equation}
respectively.
Then we construct the aggregation of the projectors over all $m$ blocks and obtain the DC estimator as
\begin{equation}
\label{eqn: aggregate_P}
\wh\bP := \frac{1}{m} \sum_{i = 1}^{m} \wh\bU_i \wh\bU_i^\top 
\in \R^{p \times p}
\text{ \ and \ }
\wh\bA_{b} := \mathrm{SVD}_{b}(\wh\bP) 
\in \R^{p \times b}, \text{ \ respectively.}
\end{equation}

Markedly, many DC or distributed PCA-type of methods \citep{Fanetal2019, GaoTsay2023} take the $r$ leading eigenvectors in their blockwise estimation step analogous to~\eqref{eqn: block_cov_vfm}.
Our investigation reveals that it is crucial to take the $b$ rather than $r$ leading eigenvectors, even when $b > r$,
so that for an uncontaminated block where $\cM_i \cap \cT = \emptyset$, the matrix $\wh\bU_i \wh\bU_i^\top$ projects onto the full column space of the block data matrix.
This allows for the exact characterization of all the eigenvalues of $\E[\wh\bU_i \wh\bU_i^\top]$ and for deriving the desired eigengap, which in turn is utilized for the estimation of the loading space by the column space of the first $r$ columns of $\wh{\bA}_b$ (see Lemma~\ref{lem: mean_block_fixed}).

To see the tail-robust property of the DC estimator, note that
retaining the full block span makes the projector $\wh\bU_i \wh\bU_i^\top$ invariant to nonzero rescaling of individual observations within the block,
and its operator norm is always bounded by one.
Therefore, anomalous observations cannot increase the operator norm of a block contribution by arbitrarily large magnitudes.
The presence of anomalies can still change the column spaces of the affected block data matrices, but only the blocks containing replacements are affected
so their influence is controlled by the fraction of affected blocks, which is captured by the term $\overline{\cT}/m$ in Theorem~\ref{thm: DC_ano_clean} below.

We point out that the individual block spaces need not consistently estimate the loading space,
and each block estimation in~\eqref{eqn: block_cov_vfm} can be viewed as a harmless ``truncation'' that could preserve the span of the loading space, before the aggregation step in~\eqref{eqn: aggregate_P}.

\medskip
\begin{remark}[Blockwise projector construction]
\label{remark: block_size}
\begin{spacing}{1.1}
Suppose the true factor number~$r$ is known. Then it is natural to set $b=r$ and hence $\wh\bA_{r} = \mathrm{SVD}_{r}(\wh\bP)$ in~\eqref{eqn: aggregate_P}.
However, the choice of retaining all $b$ eigenvectors in~\eqref{eqn: block_cov_vfm} is deliberate and essential for the desirable properties of the DC estimator.

First, our Theorem~\ref{thm: DC_ano_clean} allows for $b \geq r$, as long as the $b$ leading eigenvectors are taken for blockwise projector construction as in~\eqref{eqn: block_cov_vfm}. This is exactly the setting in which we propose the factor number estimator in Section~\ref{subsec: star_estimation}, where we set $b = r_{\max}$ with $r_{\max}$ denoting a strict upper bound on the unknown factor number $r$. 
When $r$ is known, the choice $b = r$ is the smallest block size allowable 
and it yields the best convergence rate for the DC estimator. 

On the contrary,
if $b>r$ but only the first $r$ eigenvectors are retained, i.e.,\ $\wh\bU_i = \mathrm{SVD}_{r}(\wh\bS_i)$ in~\eqref{eqn: block_cov_vfm},
then $\wh\bU_i \wh\bU_i^\top$ need not preserve the full column span of the block data matrix,
and our proof does not directly apply to this construction.
In this case, a considerably stronger moment condition may be needed to let each $\wh\bU_i \wh\bU_i^\top$ stay sufficiently close to the true projector onto the loading space so as to apply arguments as in the proof of Theorem~2 of \cite{Fanetal2019}.
It may also require the block size $b$ to be lower bounded by some unknown constant (cf.\ Theorem~4 of \citeauthor{Fanetal2019}, \citeyear{Fanetal2019}, or Theorem~1 of \citeauthor{Lietal2024}, \citeyear{Lietal2024}).
Figure~\ref{fig: sim_block_size} numerically compares this alternative with our DC construction.
\end{spacing}
\end{remark}

\begin{remark}[Connection to spatial sign covariance matrices]
\label{remark: spatial_sign_cov}
\begin{spacing}{1.1}
We note a connection between the DC estimator and robust principal component methods based on spatial sign covariance matrices.
In particular, when DC is applied with a single observation in each block such that $b=1$,
the average of all block projectors coincides with the spatial sign covariance matrix, 
$T^{-1} \sum_{t=1}^T \bx_t\bx_t^\top /\|\bx_t\|^2$.
Principal component analysis based on spatial sign covariance matrices is studied by \cite{Taskinenetal2012},
and their population eigenvalue structure under elliptical distributions is investigated by \cite{Durreetal2016}.
The correspondence between population eigenspaces does not, by itself, ensure accurate estimation when the relevant population eigengap becomes small.
For instance, under an increasingly dominant factor, the choice of $b=1$ may lead to a poor estimator in the presence of anomalies, 
whereas our DC construction with $b = r$ remains consistent under the same model even when the anomalies are present.
We show an example in the supplementary material (Section~\ref{app:rem:spatial}) that the DC estimator using singleton blocks may fail to consistently estimate the loading space
unless, for instance, stronger conditions on the factor strengths are imposed.
\end{spacing}
\end{remark}

\subsection{Split-sample tail- and anomaly-robust factor number estimator}
\label{subsec: star_estimation}

We next present the STAR method for estimating the unknown factor number $r$.
Let $r_{\max} \in [2, \min(p, \floor{T/2})]$ be an integer such that $r < r_{\max}$, i.e.,\ $r_{\max} -1$ is an upper bound on the candidate factor number.
For simplicity, we only consider the regime where
the loading matrix satisfies
\begin{equation}
\label{ass: loading_star}
c_\lambda \bI_r \preceq \bA^\top \bA /p \preceq C_\lambda \bI_r ,
\end{equation}
for some absolute constants $0 < c_\lambda \leq C_\lambda < \infty$.
Condition~\eqref{ass: loading_star} is not required for the DC loading space guarantee in Theorem~\ref{thm: DC_ano_clean} below,
and can be relaxed at the cost of more complicated arguments.
Under this, Assumption~\ref{ass: loading} holds with $\delta = 1$. 
Then, STAR proceeds as below.


\paragraph{Stage 1 (Sample splitting and DC).}
Let $\cI_1 := \{1, \dots, \lfloor T/2\rfloor\}$ and $\cI_2 := \{\lfloor T/2\rfloor + 1, \dots, T\}$ be a deterministic partition of $[T]$.
Apply the DC construction in~\eqref{eqn: block_cov_vfm} and \eqref{eqn: aggregate_P} to $\{\bx_t\}_{t\in \cI_1}$ with $b = r_{\max}$.
Denote the resulting estimator by $\wh\bA_{r_{\max}}^{(1)} \in \R^{p \times r_{\max}}$.

\paragraph{Stage 2 (projected coordinates).}
For every $t\in \cI_2$, we construct a low-dimensional projection score such that
\begin{equation}
\label{eqn: projected_coords}
\wh\bxi_t := (\wh\bA_{r_{\max}}^{(1)})^\top \bx_t \in \R^{r_{\max}} .
\end{equation}
To see why such construction can be useful, 
define $\bA_{r_{\max}} := \mathrm{SVD}_{r_{\max}}(\bSigma_x)$,
and the oracle projected vector $\bxi_t := \bA_{r_{\max}}^\top \bx_t$.
Then for a clean observation with $t\notin \cT$, note that $\bxi_t = \zeta_t \bA_{r_{\max}}^\top \bM \bv_t$ by \eqref{eqn: xt_ellip_form},
and its scatter matrix is
\[
\bA_{r_{\max}}^\top \bSigma_x \bA_{r_{\max}}
=
\diag\big( \lambda_1(\bSigma_x), \dots, \lambda_{r_{\max}}(\bSigma_x) \big)
=: \bLambda_{r_{\max}} .
\]
In particular, by Assumptions~\ref{ass: loading}, \ref{ass: noise}, 
and \eqref{ass: loading_star}, we have
\begin{equation}
\label{eqn: proj_variance}
\lambda_j(\bSigma_x) \asymp
\begin{cases}
    p, & j \in [r] , \\
    1, & r < j \leq r_{\max} .
\end{cases}
\end{equation}
Markedly, such separation between the eigenvalues of the scatter matrix of $\bxi_t$ holds even when the projected coordinates have no finite second moments, and its sample counterpart is exploited in the next stage.

\paragraph{Stage 3 (median-of-ratios).}
For each $t\in \cI_2$ and coordinate $j\in [r_{\max} -1]$, 
define the ratio
\begin{equation}
\label{eqn: ratio_def}
\wh\rho_{t,j} := \frac{\wh\xi_{t,j}^{2} }{\wh\xi_{t,j+1}^{2} } .
\end{equation}
By convention, when $\wh\xi_{t,j+1}^2 = 0$, we set $\wh\rho_{t,j} = \infty$ if $\wh\xi_{t,j}^{2} > 0$, and $\wh\rho_{t,j} = 1$ otherwise.
For the oracle projected coordinates, we can show that
$\mathrm{median}_t\big( \xi_{t,j}^2 /\xi_{t,j+1}^2 \big) = \lambda_j(\bSigma_x) /\lambda_{j+1}(\bSigma_x)$ for $j \in [r_{\max} - 1]$, which follows from the ellipticity of the data (see Lemma~\ref{lem: median_ratio_quadratic}).
From~\eqref{eqn: proj_variance}, these population medians are of order~$p$ at $j = r$ and of order~1 at all other candidate indices.
To account for potential contamination in the data, we compute the sample median of these ratios across the observations, as
\begin{equation}
\label{eqn: median_ratio}
\wh\rho_j := \mathrm{median}_{t \in \cI_2} \, \wh\rho_{t,j} ,
\qquad j \in [r_{\max} - 1] ,
\end{equation}
where we use the convention that this sample median is the $\lceil |\cI_2| / 2 \rceil$th smallest ratio.
The sample median limits the effect of anomalous projected ratios,
and the STAR estimator is set as the index at which this median ratio is maximized:
\begin{equation}
\label{eqn: r_hat_vfm}
\wh{r} := \argmax_{j \in [r_{\max} - 1]} \wh{\rho}_j .
\end{equation}
If the maximum is attained at more than one index, we take the smallest such index.

\medskip
\begin{remark}[Computational complexity]
\begin{spacing}{1.1}
For the implementation of the DC estimator, note that the local matrices $\wh\bS_i$ in~\eqref{eqn: block_cov_vfm} need not be formed explicitly.
For a block data matrix of size~$p \times b$, a thin singular value decomposition gives the required $b$ left singular vectors in $\cO(p b^2)$ operations.
Hence the total cost of the steps in~\eqref{eqn: block_cov_vfm} over all $m = \floor{T/b}$ blocks is $\cO(T p b)$.
In~\eqref{eqn: aggregate_P},
constructing $\wh\bP$ costs $\cO(T p^2)$ operations,
and the final eigendecomposition costs $\cO(p^3)$ operations.
Altogether, the computational complexity of our DC estimation is $\cO(Tpb + Tp^2 + p^3)$ which, when $b=r$ is constant, is the same as that of conventional PCA estimators which costs $\cO(Tp^2 + p^3)$.

For STAR, the DC estimation cost is of the same order as above.
Computation of the projection scores costs $\cO(T p r_{\max})$ operations,
and the subsequent ratios and medians can be computed in $\cO(T r_{\max})$ operations.
Hence the computational complexity of STAR is the same as that of the full-sample DC estimation.
\end{spacing}
\end{remark}

\section{Theoretical results}\label{sec: theorem}

For two matrices $\bA_1$ and $\bA_2$ with full column rank and the same number of rows, we measure their distance by $d(\bA_1, \bA_2) := \big\| \bA_1 (\bA_1^\top \bA_1)^{-1} \bA_1^\top - \bA_2 (\bA_2^\top \bA_2)^{-1} \bA_2^\top \big\|_F$, which can be seen as the distance between the projection spaces of $\bA_1$ and $\bA_2$.


\begin{theorem}[Statistical errors of DC]
\label{thm: DC_ano_clean}
Let Assumptions~\ref{ass: elliptical}, \ref{ass: loading}, \ref{ass: noise}, 
and \ref{ass: contamination} hold.
Recall that $m= \floor{T/b}$, and
suppose that $r \leq b$
and $b = o(p^{\min\{\delta,1\}})$. 
Let $\wh\bA_{b,1:r}$ denote the first $r$ columns of $\wh\bA_{b}$ in~\eqref{eqn: aggregate_P}.
Then there exists a constant $C > 0$
such that for all sufficiently large $p$ and every $x \in (0,1)$,
\[
\P\left( d(\wh\bA_{b,1:r}, \bA) >
C \Bigg[
p^{\max\{1 - \delta,0\}}
\left\{
\min\left(
\sqrt{\frac{b \log(4/x)}{m}} ,
\sqrt{\frac{\log(4p/x)}{m}}
\right)
+ \frac{\overline{\cT}}{m}
\right\}
+ p^{-\delta}
\Bigg] \right)
\leq x .
\]
In particular,
for any fixed $b \geq r$,
as $\min(m,p) \to \infty$, we have
\[
d(\wh\bA_{b,1:r}, \bA) =
\cO_P\left(
p^{\max\{1 - \delta,0\}}
\left\{ \frac{1}{\sqrt{m}} + \frac{\overline{\cT}}{m} \right\}
+ p^{-\delta}
\right) .
\]
Note that when $b = r$, these results directly apply to $\wh\bA_{r} := \wh\bA_{b,1:r}$.
\end{theorem}

In Theorem~\ref{thm: DC_ano_clean}, 
the $\overline{\cT}/m$ term accounts for the cost of anomalies,
while a bias-variance type of decomposition is used for analyzing the DC estimator applied to $\{\bx^0_t\}_{t \in [T]}$ as in \cite{Fanetal2019}, although here we do not require the existence of any positive moment.
Notably, if all factors are pervasive and $\delta=1$, no anomalies are present, and $b$ is fixed,
the convergence rate of DC is comparable to that of the conventional PCA estimator attainable under the assumption that the data have four or more moments, which is $\cO_P(\sqrt{1/T} + 1/p)$  \citep[e.g.,][]{BaiNg2002}.

In the scenario where $b = r$,
the rate of convergence in Theorem~\ref{thm: DC_ano_clean} becomes
$\cO_P\big( p^{1-\delta} /\sqrt{T} + p^{1-\delta}\overline{\cT} /T + p^{-\delta} \big)$ for $0 < \delta < 1$,
and $\cO_P\big( 1/\sqrt{T} + \overline{\cT}/T + p^{-\delta} \big)$ for $\delta \geq 1$.
Hence for any $\delta > 0$, a sufficient condition for the consistency of the DC estimator is
\begin{equation}
\label{eqn: DC_general_consistency}
\max\left\{ \frac{1}{p} , \;\;
\frac{p^{2 \max\{1 - \delta,0\}}}{T} , \;\;
\frac{p^{\max\{1 - \delta,0\}} \overline{\cT}}{T} \right\}
= o(1) .
\end{equation}

Next, we make the following assumption for deriving the consistency of the STAR estimator.

\begin{assumption}
\label{ass: sample_split_independence}
Recall the partition $\cI_1$ and $\cI_2$ in Section~\ref{subsec: star_estimation}.
We assume the observations $\{\bx_t\}_{t\in \cI_1}$ and $\{\bx_t\}_{t\in \cI_2 \setminus \cT}$ are independent,
conditional on $\cT$.
\end{assumption}

Together with Assumption~\ref{ass: contamination}, Assumption~\ref{ass: sample_split_independence} ensures that the second-half clean observations retain their original distribution conditionally on the first-half data and the anomaly locations.
Note that it does not rule out the case where anomalous observations in $\cI_1$ depend on the clean observations in $\cI_1$,
while imposing no additional restriction on the anomalous observations in $\cI_2$.

\begin{theorem}[Consistency of STAR]
\label{thm: r_hat_consistency}
Let Assumptions~\ref{ass: elliptical}, \ref{ass: loading}, \ref{ass: noise}, 
\ref{ass: contamination},
\ref{ass: sample_split_independence},
and \eqref{ass: loading_star} hold.
Suppose that $r < r_{\max}$ and
$b = r_{\max} = o(p)$.
Recalling that $m = \floor{T/b}$, we further assume that 
\begin{equation}
\label{eqn: thm_r_contam_condition}
\sqrt{p} \log(p)
+ \overline{\cT} \, p^{1/4}
= o(m).
\end{equation}
Then as $\min(m,p) \to \infty$,
the STAR estimator in~\eqref{eqn: r_hat_vfm} is consistent such that $\P(\wh{r} = r) \to 1$.
\end{theorem}

Theorem~\ref{thm: r_hat_consistency} establishes the consistency of the proposed STAR estimator.
Condition~\eqref{eqn: thm_r_contam_condition} ensures separation of the conditional population medians of the clean score ratios and controls the effect of sampling fluctuations and anomalies on their empirical medians.
More concretely, for each $j \in [r_{\max}-1]$, consider the median of the conditional distribution of $\wh\rho_{t,j}$ for $t \in \cI_2 \setminus \cT$, given the sigma-algebra $\sigma(\cT, \{\bx_t\}_{t\in\cI_1})$.
Then with probability approaching one, these conditional medians are uniquely maximized at $j = r$.

\section{Estimation of the factor structure}
\label{sec: factor_structure}

When the factor number is unknown, upon the consistent estimation of $r$ by the STAR estimator $\wh r$ (Theorem~\ref{thm: r_hat_consistency}), we can perform the DC estimation again to all observations, now with block size $b = \wh{r}$.
Denoting the resulting estimator by $\wh\bA_{\wh{r}}$, the following guarantee follows immediately.

\begin{corollary}[Convergence rate of DC combined with STAR]
\label{cor: dc_star}
Let all conditions in Theorem~\ref{thm: r_hat_consistency} hold.
Then as $\min(T,p) \to \infty$, we have
\[
d(\wh\bA_{\wh{r}}, \bA) =
\cO_P\left( 
\frac{1}{\sqrt{T}}
+ \frac{1}{p}
+ \frac{\overline{\cT}}{T}
\right) .
\]
\end{corollary}

Next, we discuss the estimation of the latent core factor or the factor-driven common component under~\eqref{eqn: vfm}. 
For this, we suppose that we directly observe $\{\bx^0_t\}_{t \in [T]}$ or, equivalently, that $\cT = \emptyset$.
Then with the DC estimator $\wh\bA_{\wh{r}}$, for every $t\in[T]$, we can subsequently estimate the core factor $\bff_t$ and thus the common component $\bchi_t = \bA \bff_t$, respectively, as
\begin{equation}
\label{eqn: dc_fac_common_est}
\wh\bf_t := \frac{1}{\sqrt{p}} \wh\bA_{\wh{r}}^\top \bx_t,
\qquad
\wh\bchi_t := \wh\bA_{\wh{r}} \wh\bA_{\wh{r}}^\top \bx_t ,
\end{equation}
where $\wh\bchi_t$ is essentially the orthogonal projection of $\bx_t$ onto the estimated loading space.

\begin{corollary}[Convergence rates for factor and common component estimation]
\label{cor: dc_fac_common}
Suppose that all the conditions in Theorem~\ref{thm: r_hat_consistency} hold.
Assume $\cT = \emptyset$, 
and $\zeta/\sqrt{p} = \cO_P(1)$ as $\min(T,p) \to \infty$.
Define the matrix $\wh\bH := \wh\bA_{\wh{r}}^\top \bA /\sqrt{p}$.
Then as $\min(T,p) \to\infty$,
the singular values of $\wh\bH$ are bounded away from zero and infinity with probability approaching one, 
and for every $t\in[T]$,
\begin{align*}
    \big\| \wh\bf_t - \wh\bH \bf_t \big\|^2 &= 
    \cO_P\left(
    \frac{1}{T}
    +\frac{1}{p} 
    \right) , 
    \qquad
    \frac{1}{p} \big\| \wh\bchi_t - \bchi_t \big\|^2 =
    \cO_P\left(
    \frac{1}{T}
    +\frac{1}{p} 
    \right) .
\end{align*}
\end{corollary}

The additional condition that $\zeta/\sqrt{p} = \cO_P(1)$ is imposed to control the stochastic scale of individual observations, see also Remark~\ref{remark: example_scale_wo_moment}.
A similar condition is found in Assumption~A of \cite{HeKongYuZhang2022}, 
and a stronger version, viz.\ $\E[\zeta^2/p] = 1$, is found in Assumption~4.1 of \cite{FanLiuWang2018}.

\medskip
\begin{remark}\label{remark: example_scale_wo_moment}
\begin{spacing}{1.1}
The condition $\zeta_t/\sqrt{p} = \cO_P(1)$ does not require any finite positive moment.
For example, let $Z_t$ be independent copies of a random variable with density $f_Z(z) = z^{-2} \mathbbm{1}_{\{z \geq 1\}}$, and be independent of the direction vectors in Assumption~\ref{ass: elliptical}, 
and set $\zeta_t = \sqrt{p} \exp(Z_t)$.
Then $\zeta_t$ is finite and strictly positive almost surely as required in Assumption~\ref{ass: elliptical}.
As $M \to \infty$, we have
\[
\sup_{p \geq 1}
\P\left( \frac{\zeta_t}{\sqrt{p}} > M \right)
=
\P(Z_t > \log(M))
=
\frac{1}{\log(M)}
\to 0 ,
\]
which shows $\zeta_t/\sqrt{p} = \cO_P(1)$.
However, for every $q > 0$, we have
\[
\E\left[ \left( \frac{\zeta_t}{\sqrt{p}} \right)^q \right]
=
\int_1^\infty
\frac{\exp(qz)}{z^2} \,\mathrm{d}z
= \infty .
\]
\end{spacing}
\end{remark}

\section{Numerical experiments}
\label{sec: numerical}

\subsection{Simulation}
\label{subsec: sim}

We investigate the finite-sample performance of the proposed estimators.
Our data generating mechanism is motivated by that of \cite{HeKongYuZhang2022}.
In particular, we include a continuous parameter which controls the departure from ellipticity together with an autoregressive coefficient to both factors and noise, so as to assess the estimation performance when the underlying model departs from our assumptions.

\paragraph{Data generating mechanism.}
Throughout, we fix $r = 2$; other specifications examined in our preliminary experiments result in similar patterns.
For factor number estimation by STAR, we set $r_{\max} = 5$ and hence $b = r_{\max}$.
For loading space estimation by DC with the true factor number given, we set $b = r = 2$.
We detail the data generation below.
First, independently across $t$, generate $\bg_t\sim \cN(\zero,\bI_{r+p})$ and $\bs_t\in \R^{r+p}$ whose coordinates are independent copies of $E-1$, where $E\sim \mathrm{Exp}(1)$.
The vectors $\bg_t$ and $\bs_t$ are independent.
For some $a\in [0,1]$, define $\bq_t := \sqrt{1-a}\, \bg_t + \sqrt{a}\, \bs_t$
and $\bw_t := R_t \bq_t$,
where $R_t$ is independent of these vectors and independent across time, with
$R_t := \sqrt{(\nu-2) /V_t}$ for $\nu>2$,
$R_t = \sqrt{\nu/V_t}$ for $\nu \in (0,2]$,
and $V_t \sim \chi_\nu^2$.
The variables $V_t$ are independent across $t$ and independent of $\{\bg_t, \bs_t\}_{t\in [T]}$.

For $\nu > 2$, this construction gives $\E[\bw_t] = \zero$ and $\E[\bw_t \bw_t^\top] = \bI_{r+p}$ for every $a$.
When $a=0$, the vector $\bw_t$ is spherical with either a covariance-normalized multivariate $t_\nu$ distribution if $\nu>2$,
or a standard multivariate $t_\nu$ distribution if $\nu\in (0,2]$.
Note that the latter does not have finite second moments, and even its mean is undefined when $\nu \leq 1$.
Setting $a>0$ introduces asymmetry and violates the setup of elliptical factor models.
With the above,
for $|\phi|<1$,
we generate the AR(1) process $\bh_t = \phi\bh_{t-1} + \sqrt{1-\phi^2} \, \bw_t$,
where we set $\bh_1=\bw_1$ in implementation and only retain observations after a short period of initialization.

The loading matrix is set as $\bA = \sqrt{p} \,\bQ$, where $\bQ \in \R^{p \times r}$ is the half-orthogonal matrix from QR decomposition of a $p\times r$ matrix with $\cN(0,1)$ entries independent of each other and independent of the innovations.
Also, with $\bSigma_\epsilon = (0.3^{|j-k|})_{j,k \in [p]}$, let $\bSigma_\epsilon = \bA_\epsilon \bA_\epsilon^\top$ be its Cholesky decomposition with $\bA_\epsilon$ lower triangular.
Then we partition $\bh_t=(\bh_{f,t}^\top, \bh_{\epsilon,t}^\top)^\top$ with $\bh_{f,t}\in\R^r$, and generate the clean observations as
\[
\bf_t = \bh_{f,t},
\qquad
\bepsilon_t = \bA_\epsilon \bh_{\epsilon,t},
\qquad
\bx^0_t = \bA\bf_t + \bepsilon_t .
\]
Note that for $\nu>2$ and all choices of $a$ and $\phi$, we have $\E[\bf_t\bf_t^\top] = \bI_r$,
$\E[\bepsilon_t\bepsilon_t^\top] = \bSigma_\epsilon$,
and $\E[\bf_t\bepsilon_t^\top] = \zero$.
The lag-$h$ covariance of the clean observation $\bx^0_t$ is then $\phi^h(\bA\bA^\top + \bSigma_\epsilon)$ for $h\geq 0$.
For $a = 0$ and $\phi = 0$, all values $\nu > 0$ yield clean observations satisfying the moment-free elliptical formulation.

To introduce anomalies, we select $\cT \subset [T]$ of size $\lceil\varepsilon T\rceil$ uniformly and independently of the clean process.
For each $t\in\cT$, independently generate
$\bz^0_t\sim N(\zero,\bI_p)$, set $\bz_t = 100 \sqrt{p}\, \bz^0_t /\|\bz^0_t\|$ and generate $\bx_t$ as in~\eqref{eq:obs},
where $\varepsilon \in [0,1/2)$ controls the fraction of anomalies.

\begin{figure}[ht!]
\centering
\includegraphics[width=.85\linewidth]{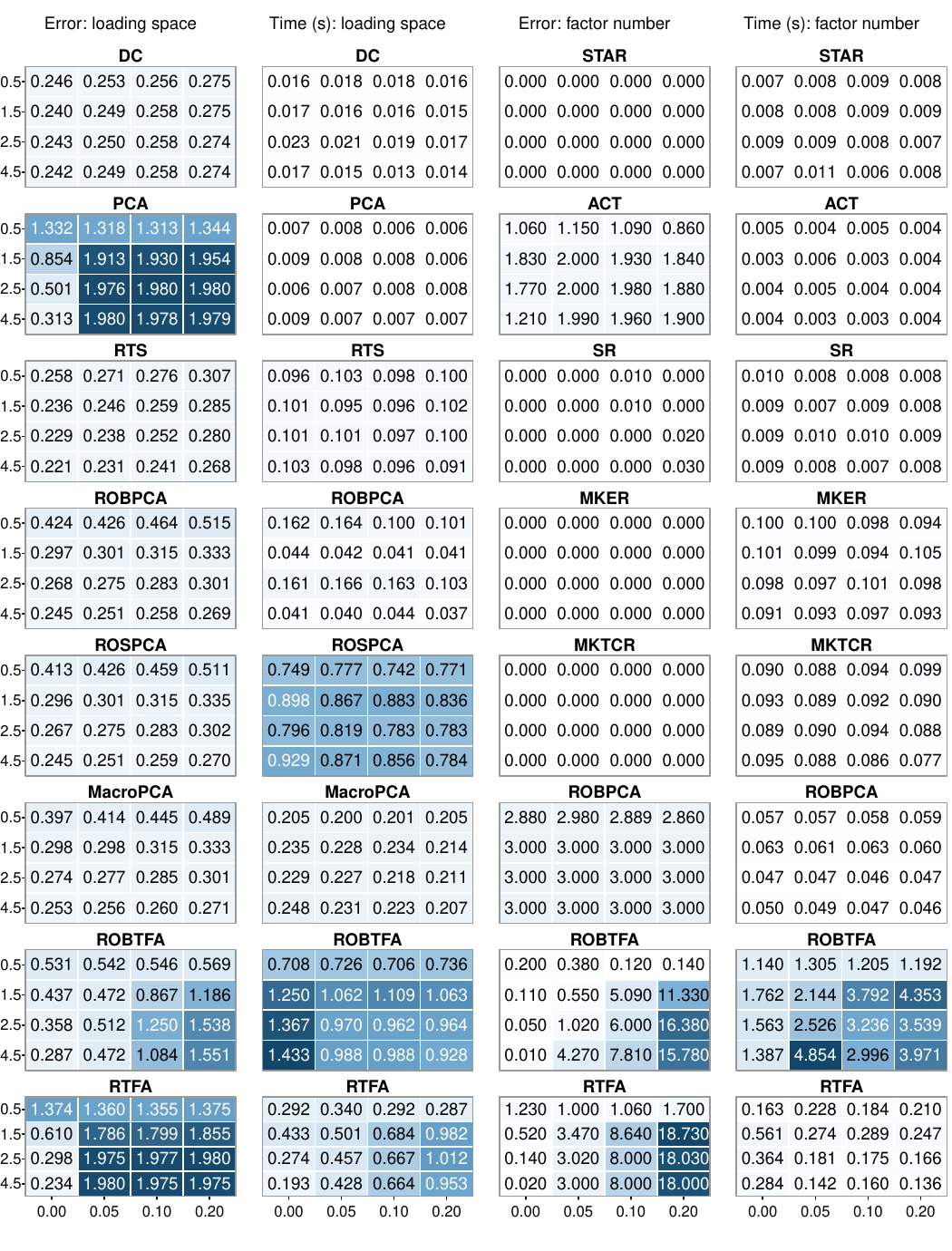}
{\setstretch{1}
\caption{
Performance of loading space and factor number estimation for $(T,p)=(100,100)$.
Within each heatmap, rows correspond to
$\nu \in \{ 0.5, 1.5, 2.5, 4.5 \}$
and columns correspond to
$\varepsilon \in \{ 0, 0.05, 0.1, 0.2 \}$.
Darker shading indicates larger values within each column.
}
\label{fig: sim_nu_epsilon}}
\end{figure}

\paragraph{Competing methods.}
For each competing method below, any tuning parameters are specified as suggested in the corresponding work.
For loading space estimation, we provide the true factor number $r$ to each method, and compare DC with:
(i)~the plain principal component estimator \citep[e.g.][]{BaiNg2002, Bai2003}, denoted by \texttt{PCA};
(ii)~the robust two-step method of \cite{HeKongYuZhang2022}, denoted by \texttt{RTS};
(iii)~\texttt{ROBPCA} of \cite{Hubertetal2005}, which combines projection pursuit with robust scatter matrix estimation;
(iv)~\texttt{ROSPCA} of \cite{Hubertetal2016}, which incorporates sparse loading estimation into a robust principal component procedure;
(v)~\texttt{MacroPCA} of \cite{Hubertetal2019}, which accommodates cellwise outliers;
(vi)~\texttt{ROBTFA} of \cite{Barigozzietal2026_robust},
based on entrywise truncation;
and (vii)~\texttt{RTFA} of \cite{Barigozzietal2023_RTFA}, based on a Huber loss applied to the normalized residual norm iteratively. 

For estimating the factor numbers, we compare STAR with:
(i)~the adjusted correlation thresholding estimator of \cite{FanGuoZheng2022}, denoted by \texttt{ACT};
(ii)~the Spearman correlation matrix based estimator of \cite{Qiuetal2025}, denoted by \texttt{SR};
(iii)--(iv)~the multivariate Kendall's tau eigenvalue ratio and transformed contribution ratio estimators of \cite{Yuetal2019}, denoted by \texttt{MKER} and \texttt{MKTCR}, respectively;
(v)~the automatic component selection procedure associated with \texttt{ROBPCA} of \cite{Hubertetal2005};
(vi)~the data truncation-based rank estimator of \cite{Barigozzietal2026_robust}, denoted by \texttt{ROBTFA};
and (vii)~the Huber loss-based eigenvalue ratio estimator of \cite{Barigozzietal2023_RTFA}, denoted by \texttt{RTFA}.

\paragraph{Settings and results.}
To measure the performance of different estimators,
we provide the true factor number~$r$ to each method and report the distance $d(\cdot, \cdot)$ introduced earlier for loading space estimation,
and we report the absolute error $|\wh{r} - r|$ for factor number estimation.
Each error is averaged over $100$ Monte Carlo replications, and reported with the average computation time in seconds for the corresponding estimation task.
To align with the setup in this paper, we first consider $a = 0$ and $\phi = 0$, so that the clean observations are independent across time and follow an elliptical distribution.
We vary the degrees of freedom over $\nu \in \{0.5, 1.5, 2.5, 4.5\}$ and the fraction of anomalies over $\varepsilon \in \{0, 0.05, 0.1, 0.2\}$.
The former grid includes the distributions without finite first, second, and fourth moments when $\nu=0.5$, $1.5$, and $2.5$, respectively.

We present the results for $(T,p)=(100,100)$ in Figure~\ref{fig: sim_nu_epsilon}; 
additional experiments are reported in the supplement (Section~\ref{add:sim}) and exhibit similar patterns.
DC maintains small errors which increase only moderately with $\varepsilon$ and vary little over $\nu$.
Notably, its errors are either the smallest or comparable to most competing methods such as
\texttt{RTS}, \texttt{ROBPCA}, \texttt{ROSPCA}, and \texttt{MacroPCA}.
While it is unsurprising that \texttt{PCA} is sensitive to both heavy tails and anomalies,
\texttt{RTFA} also deteriorates in the presence of anomalies and \texttt{ROBTFA} exhibits a gradual increase in error as the fraction of anomalies grows.
For factor number estimation, STAR, \texttt{MKER}, and \texttt{MKTCR} recover the true factor number in every replication across the settings considered.
\texttt{SR} also performs well, with small errors appearing under anomalies.
The behavior of \texttt{ROBTFA} and \texttt{RTFA} are akin to their loading space estimation.
The performance of \texttt{ACT} is as expected since it requires more than $6$th moment (cf.\ Assumption~(b) of \citeauthor{FanGuoZheng2022}, \citeyear{FanGuoZheng2022}).
Method~\texttt{ROBPCA} (in particular, its rank selection procedure here) performs poorly, as it is only proposed aside its subspace estimation.

A salient feature of DC is its computation time, which is the lowest in almost all settings, 
among all competing methods
excluding those designed for light-tailed data (i.e.,\ \texttt{PCA} for loading space estimation and \texttt{ACT} for factor number estimation).
Particularly for loading space estimation, while DC is slower but remains comparable to \texttt{PCA}, it is significantly faster than all other methods.

\paragraph{Sensitivity to the block size.}
We examine the role of the block size and the local retained rank, as discussed in Remark~\ref{remark: block_size}.
We fix $(T,p) = (100,100)$, $r = 2$,
$a = \phi = 0$, $\varepsilon = 0$,
and consider $\nu \in \{0.5, 1.5, 2.5, 4.5\}$.
For each configuration, we vary the block size over $b\in \{2, 4, 8, 16, 32, 64, 100\}$.
We compare the proposed DC construction, which retains the leading $b$ eigenvectors of each $\wh\bS_i$ in~\eqref{eqn: block_cov_vfm}, 
with an alternative that retains only the leading $r$ eigenvectors, irrespective of $b$.
Both procedures average the resulting local projectors and retain the leading $r$ eigenvectors of the aggregate.
The two procedures coincide when $b = r$,
and the alternative would reduce to \texttt{PCA} when $b = T$.

\begin{figure}[ht!]
\centering
\includegraphics[width=.75\linewidth]{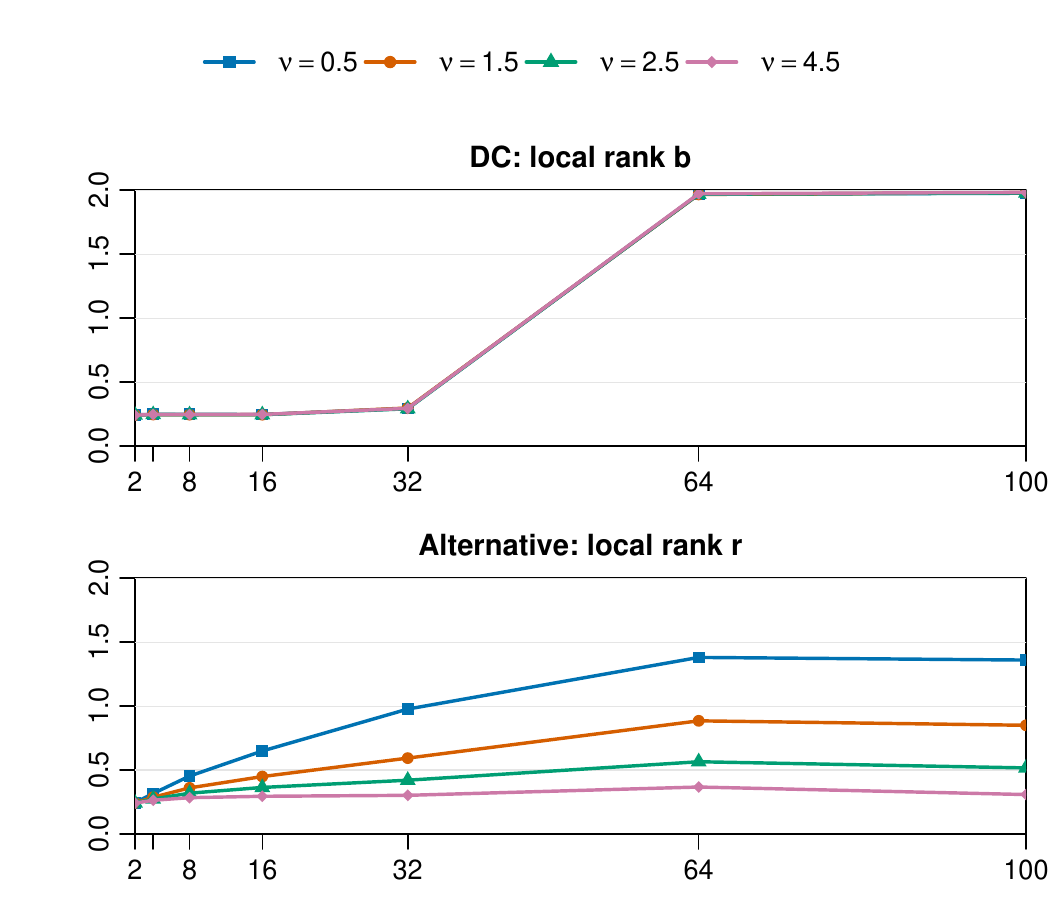}
{\setstretch{1}
\caption{
Sensitivity to the block size for $(T,p) = (100,100)$, $r = 2$, $a = \phi = 0$ and $\varepsilon = 0$.
The upper panel retains $b$ eigenvectors locally as in our DC construction;
the lower panel retains $r$ eigenvectors locally.
The horizontal axis is over the block size $b$, and the vertical axis gives the average loading space error $d(\cdot,\bA)$ over $100$ replications.
Different curves correspond to
$\nu \in \{0.5, 1.5, 2.5, 4.5\}$.
}
\label{fig: sim_block_size}}
\end{figure}

Figure~\ref{fig: sim_block_size} illustrates the distinction between the two local constructions.
For our proposed DC construction, the curves nearly coincide across $\nu$, and the errors remain small for $b \in \{2, 4, 8, 16, 32\}$.
Such insensitivity to the tail parameter pinpoints why we should retain the full block span.
The sharp increase in the upper panel at $b = 64$ and $b = 100$ is due to the number of blocks $m = 1$, so that the aggregate is a single rank-$b$ projector whose first $b$ eigenvalues are all one.
Consequently, its leading rank-$r$ subspace is not uniquely determined and chosen randomly.

For the alternative construction, retaining only $r$ directions when $b > r$ no longer preserves the full block span.
Its errors increase markedly with $b$ under the heavier-tailed configurations, particularly when $\nu = 0.5$ and $1.5$,
whereas the deterioration is much milder when $\nu = 4.5$.
At $b = T$, this alternative coincides with \texttt{PCA}, and thus gives performance as expected.
Comparing the two panels thus supports the proposed blockwise projector construction for the DC estimator.

\begin{figure}[ht!]
\centering
\includegraphics[width=.82\linewidth]{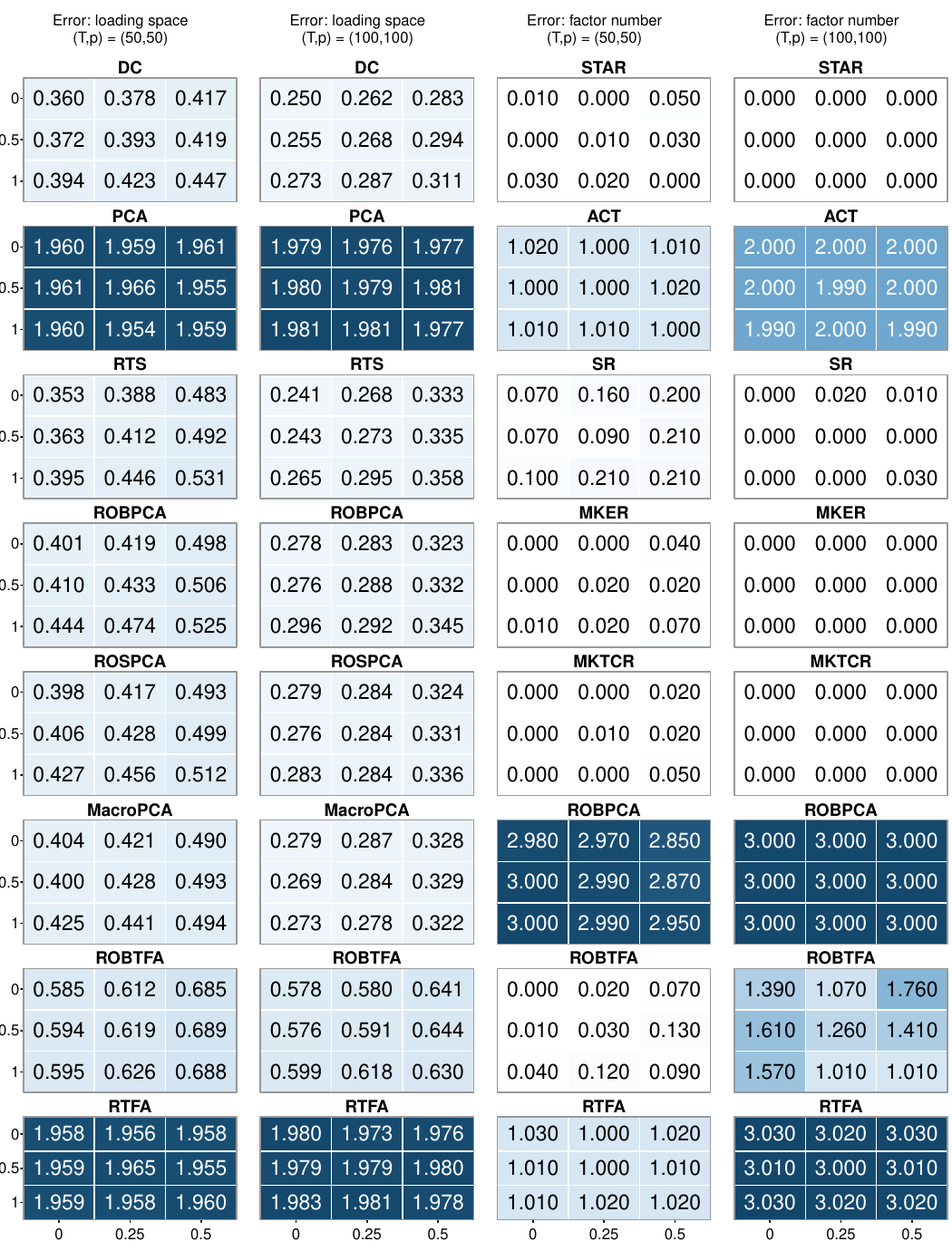}
{\setstretch{1}
\caption{
Performance of loading space and factor number estimation under non-ellipticity and temporal dependence,
with $\nu = 2.5$ and $\varepsilon = 0.05$.
Within each heatmap, rows correspond to
$a \in \{0, 0.5, 1\}$
and columns correspond to
$\phi \in \{0, 0.25, 0.5\}$.
Darker shades indicate larger values within each column.
}
\label{fig: sim_sens}}
\end{figure}

\paragraph{Model misspecification.}
We next investigate sensitivity to non-ellipticity and temporal dependence.
We fix $\nu = 2.5$ and $\varepsilon = 0.05$, 
and consider all combinations of $a\in \{0, 0.5, 1\}$, $\phi \in \{0, 0.25, 0.5\}$, and $(T,p) \in \{(50,50), (100,100)\}$,
where $a$ controls non-ellipticity and $\phi$ controls temporal dependence in the generated data.
All remaining parameters follow the preceding experiment. 
Each setting is repeated $100$ times, and the results are summarized in Figure~\ref{fig: sim_sens}.

From the figure, DC remains accurate across the parameter grids,
and increasing either $a$ or $\phi$ leads to a moderate increase in its error, while increasing $T$ and $p$ jointly reduces the error in every corresponding setting.
DC and \texttt{RTS} exhibit comparable performance when $\phi = 0$, with RTS attaining slightly smaller errors in most of these settings.
However, when $\phi \in \{0.25,0.5\}$, DC has smaller errors than \texttt{RTS} throughout the grid, suggesting that the former is less sensitive to temporal dependence in this experiment.
The errors of \texttt{PCA} and \texttt{RTFA} remain close to $2$, indicating the failure of these methods in these scenarios.

For factor number estimation, STAR remains accurate across all combinations of $a$ and $\phi$.
At $(T,p) = (50,50)$, its mean absolute rank error is at most $0.050$, and all errors are zero at $(T,p) = (100,100)$.
\texttt{MKER} and \texttt{MKTCR} exhibit similarly strong performance.
Compared to the good performance of \texttt{SR} in Figure~\ref{fig: sim_nu_epsilon}, it appears quite sensitive to the non-ellipticity and temporal dependence, which is interesting since this method does not require ellipticity \citep{Qiuetal2025}.
These results provide numerical evidence that DC and STAR remain effective under the departures from ellipticity and temporal independence considered here.

Overall, the experiments numerically show the effectiveness of DC and STAR when compared to a variety of competing methods.
Due to space constraints, we defer additional experiments on common component estimation in the supplementary material (Section~\ref{add:sim}).

\begin{figure}[ht!]
\centering
\includegraphics[width = \linewidth]{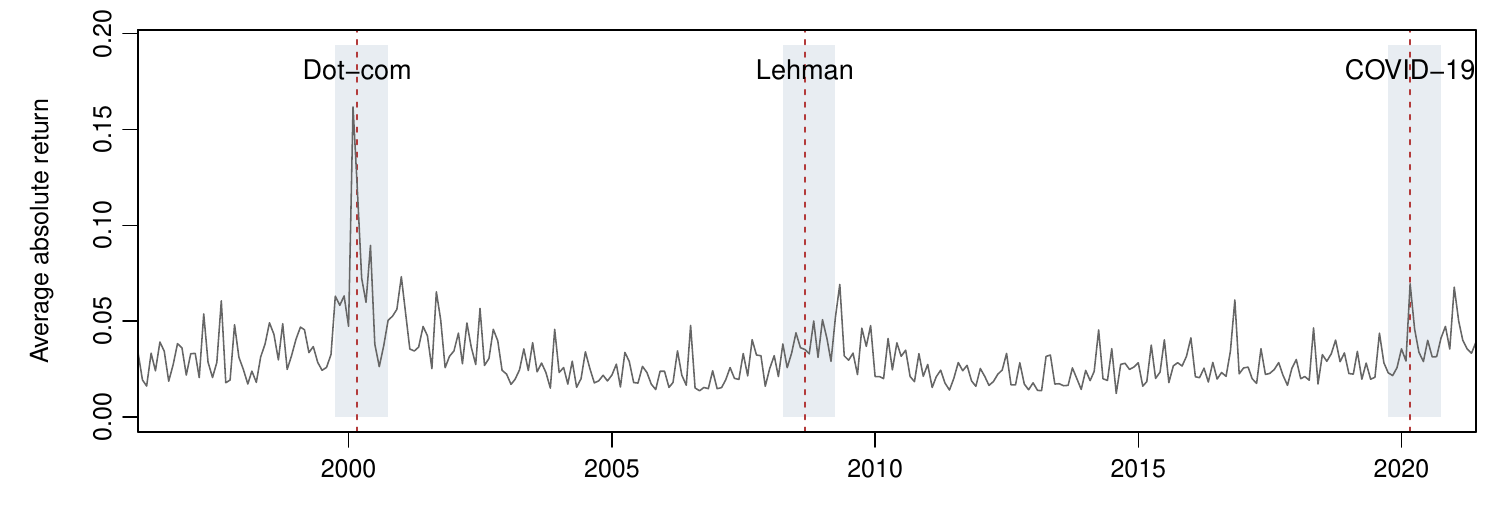}
\caption{
Average absolute return, $p^{-1} \sum_{j=1}^p |x_{t,j}|$, for the portfolio data from January 1996 to June 2021.
Shaded regions indicate the one-year windows removed in the separate comparisons.
Vertical dashed lines mark March 2000, September 2008, and March 2020.
}
\label{fig: real_returns}
\end{figure}

\subsection{Real data analysis}
\label{subsec: real_data}

In this subsection, we investigate the stability of the proposed methods using the Fama--French portfolio returns.
Particularly, we analyze monthly observations from January 1996 to June 2021 for $100$ portfolios formed by sorting stocks on size and book-to-market equity, and denote each observation by $\bx_t \in \R^p$ with $p = 100$ and $T = 306$.
See the details of the data preparation in Section~\ref{subsec: real_data_details} of the supplementary material.

\begin{figure}[ht!]
\centering
\includegraphics[width = .9\linewidth]{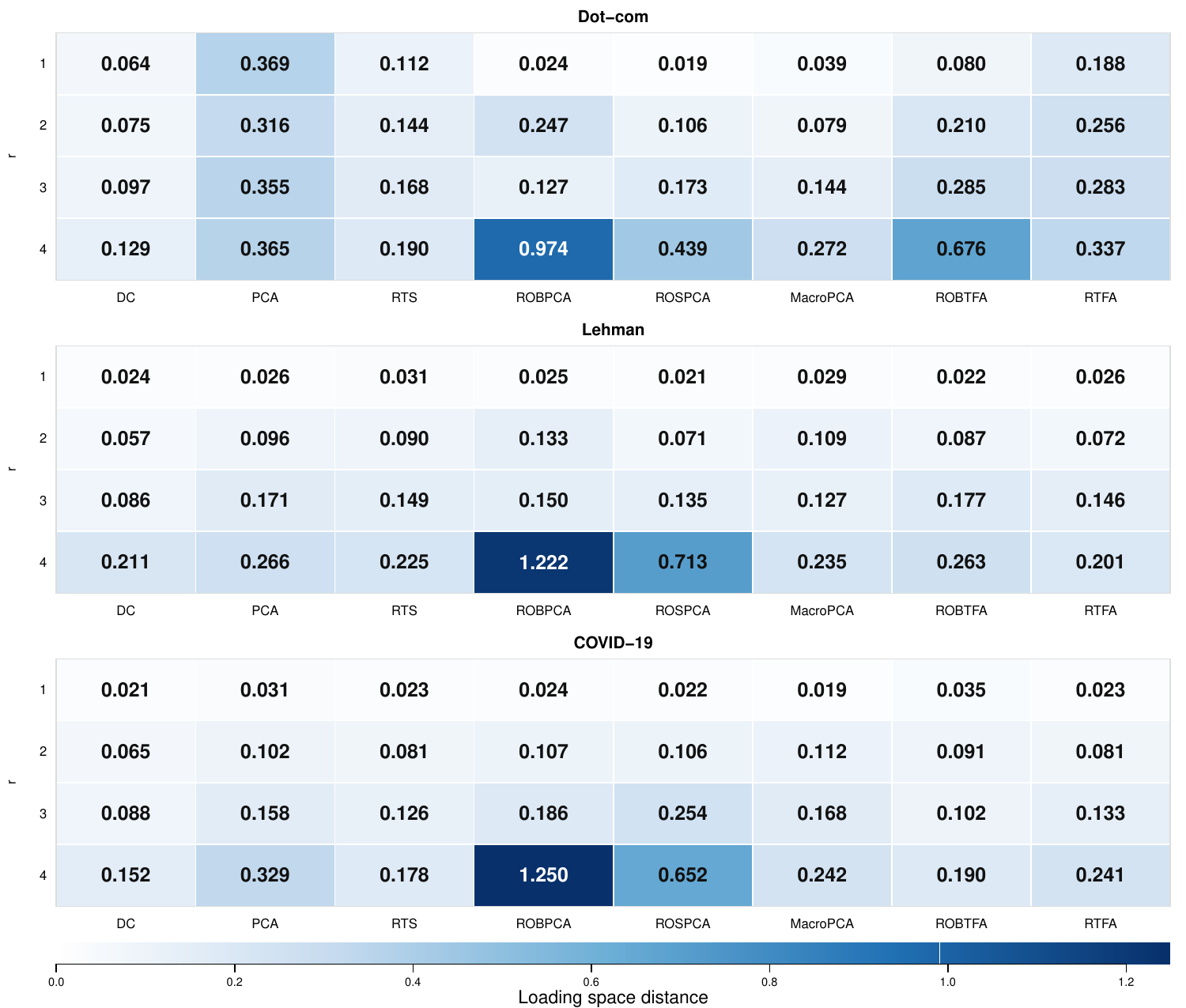}
\caption{
Loading space sensitivity to separate exclusions of the three events.
Within each panel, rows correspond to $r\in \{1, 2, 3, 4\}$ and columns correspond to the eight methods.
Darker shading indicates larger values, which represent lower stability.
}
\label{fig: real_loading_stability}
\end{figure}

We assess the sensitivity of the estimation results to the three periods associated with major market disruptions:
(i) October 1999--September 2000 for the dot-com bubble;
(ii) April 2008--March 2009 for the collapse of Lehman Brothers;
and (iii) October 2019--September 2020 for the onset of COVID-19.
For each event, we remove a one-year window and compare the estimates obtained from the remaining observations with those obtained from the complete period.
The three windows are removed one at a time, so each reduced sample contains $294$ observations.
Figure~\ref{fig: real_returns} displays the average absolute return and the selected windows,
where plot exhibits a pronounced spike around the dot-com bubble, with smaller elevations around the other two events.

First, we compare loading space estimates at each rank specified from $r\in \{1, 2, 3, 4\}$.
For every method and event, we compute the loading space distance $d(\cdot,\cdot)$ between the rank-$r$ estimate from the complete and reduced samples.
Figure~\ref{fig: real_loading_stability} shows that DC remains stable overall, with its distances ranging from $0.021$ to $0.211$, and are either comparable to competing methods or smaller than e.g.\ those of \texttt{PCA} and \texttt{RTS} in every displayed comparison.
This indicates that the loading spaces estimated by DC are less sensitive to the removal of potentially anomalous periods.

We next examine factor number estimation, using $r_{\max} = 5$ for STAR, so that its candidate factor numbers are $1,\dots,4$, and the same for competing methods.
The detailed results are relegated to the supplementary material, 
which show that all the factor number estimators retain their respective estimates when excluding each of the three events.
In particular, STAR, \texttt{SR}, \texttt{MKER}, \texttt{ROBTFA}, and \texttt{RTFA} select one factor, 
whereas \texttt{ACT}, \texttt{MKTCR}, and \texttt{ROBPCA} select more than one factors.

\begin{figure}[ht!]
\centering
\includegraphics[width = \linewidth]{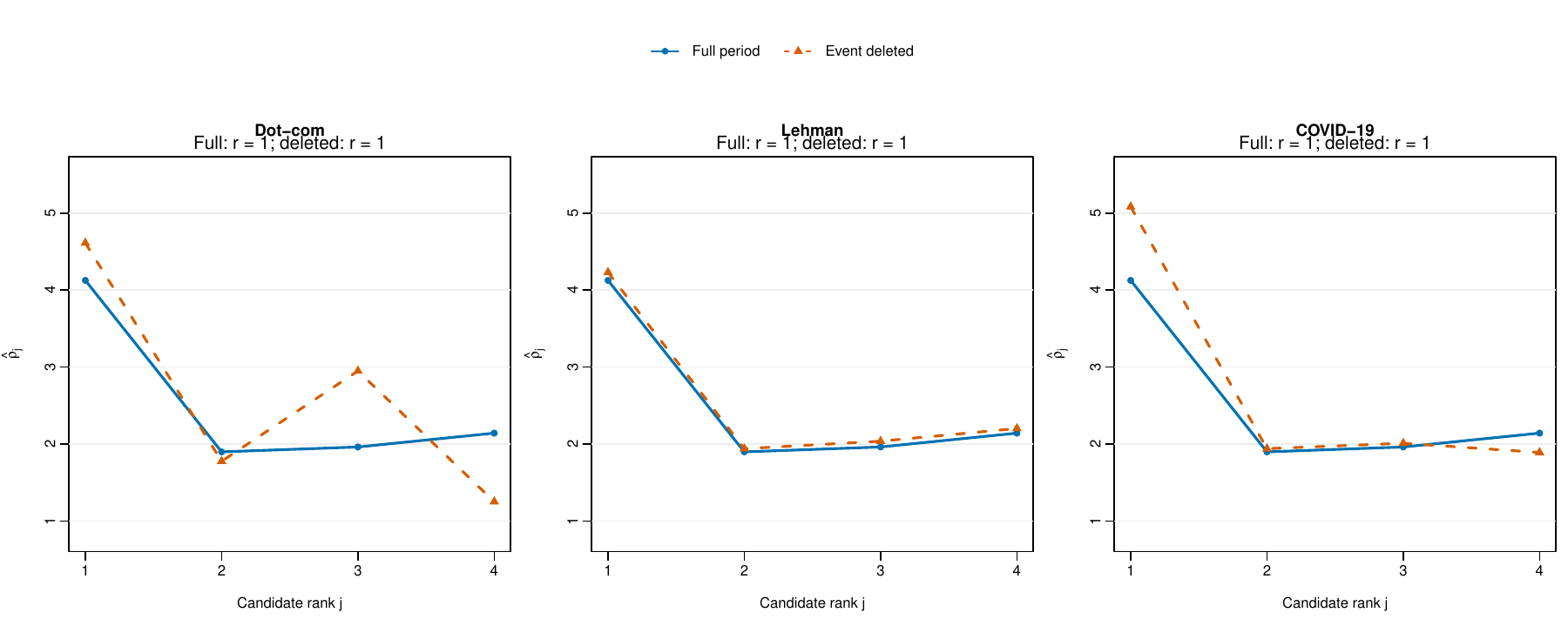}
\caption{
Median-ratio criterion of STAR with $r_{\max} = 5$.
Each panel compares the complete sample with the sample obtained by removing the indicated event window.
The vertical axis shows the median-ratio criterion $\wh\rho_j$ in~\eqref{eqn: median_ratio},
for $j\in \{1, 2, 3, 4\}$.
}
\label{fig: real_star_profiles}
\end{figure}

For STAR, we further show in Figure~\ref{fig: real_star_profiles} the corresponding ratio criterion, akin to the scree plots for the ratio-based estimators as in e.g.,\ \cite{ahn2013eigenvalue}.
The diagnostic plots display the median ratios directly, with the maximizing index giving the estimated factor number.
The considerable spikes in the first index again highlights the stability of STAR.

\section{Concluding remarks}
\label{sec: conclusion}

We have studied loading space and factor number estimation under a moment-free elliptical factor model with anomalies.
The DC estimator uses the full column spaces of small blocks before aggregation,
and we derive an explicit bound on loading space error that accounts for both sampling fluctuation and contamination.
The STAR estimator combines DC with medians of projected coordinate ratios, and is shown to consistently estimate the factor number.
The resulting loading space estimator can then be used for factor and common component estimation.

Our numerical experiments suggest that the robustness property (with respect to temporal dependence, non-ellipticity) of DC might be more ample than what we have established in this paper, which merits further investigation.
Also, we note that while elliptical distribution is inevitable in the technical arguments, our results may serve as an intermediate step in controlling the rates of performance if one can characterize how far the data distribution is from the set of all elliptical factor models.


\bibliographystyle{apalike}
\bibliography{reference.bib}


\clearpage
\appendix
\renewcommand{\thesection}{\Alph{section}}
\renewcommand{\theHsection}{\Alph{section}}
\renewcommand{\thesubsection}{\Alph{section}.\arabic{subsection}}
\renewcommand{\theHsubsection}{\Alph{section}.\arabic{subsection}}


\section{Proofs of theorems and auxiliary results}

Let the half-orthogonal matrix $\ddot\bA \in\R^{p\times r}$ represent an orthonormal basis of $\spa(\bA)$.
Hence throughout this section,
for any half-orthogonal matrix $\wh\bA$, we consider
\[
d(\wh\bA, \bA) = \big\| \wh\bA (\wh\bA^\top \wh\bA)^{-1} \wh\bA^\top - \bA (\bA^\top \bA)^{-1} \bA^\top \big\|_F 
= 
\big\| \wh\bA \wh\bA^\top - \ddot\bA \ddot\bA^\top \big\|_F 
=
d(\wh\bA, \ddot\bA) .
\]

Throughout, let $\bSigma_x = \bV_x \bLambda_x \bV_x^\top$ denote the eigendecomposition of $\bSigma_x$,
where $\bV_x$ is orthogonal 
and $\bLambda_x = \diag(\lambda_1(\bSigma_x), \dots, \lambda_p(\bSigma_x))$ has its diagonal entries in non-increasing order.
Also define $\bA_r^\ast := \mathrm{SVD}_r(\bSigma_x)$
and $\bPi_r^\ast := \bA_r^\ast (\bA_r^\ast)^\top$.
Note that by Assumptions~\ref{ass: loading} and \ref{ass: noise},
the matrix $\bPi_r^\ast$ is uniquely defined for all sufficiently large $p$.




\subsection{Proof of theorems}
\label{subsec: proof_thm}

\begin{proof}[Proof of Theorem~\ref{thm: DC_ano_clean}]
Recall the latent observations $\bx_t^0$ in~\eqref{eqn: vfm}.
For each~$i\in [m]$, let $\bP_i^0$ be the orthogonal projector onto the column space of the $i$th block $[\bx_t^0: t\in \cM_i]$,
and define $\wh\bP^0 := m^{-1} \sum_{i=1}^m \bP_i^0$
and $\bP^\ast := \E[\bP_i^0]$.
Then by Assumption~\ref{ass: elliptical} and Lemma~\ref{lem: block_projector_unique},
the matrices $\bP_i^0$, $i\in[m]$, are i.i.d.\ and have rank~$b$ almost surely.
Hence by Lemma~\ref{lem: subgauss_aggregate_P}, with probability $x/2$ in each of its two bounds, we obtain
\begin{equation}
\label{eqn: DC_ano_high_prob}
\big\| \wh\bP^0 - \bP^\ast \big\|
\leq 
C \min\left\{
\sqrt{\frac{b \log(4/x)}{m}} ,
\sqrt{\frac{\log(4p/x)}{m}}
\right\} ,
\end{equation}
with probability at least $1 - x$,
where for the second bound in Lemma~\ref{lem: subgauss_aggregate_P}, we used the fact that
\begin{align*}
    \sqrt{\frac{b}{m}} 
    + \sqrt{\frac{b \log(2/x)}{m}}
    &=
    \sqrt{\frac{b}{m}}
    \left\{ 1 + \sqrt{\log(2/x)} \right\}
    \leq
    2 \sqrt{\frac{b \log(4/x)}{m}} .
\end{align*}

Denote by $\cB:= \{i\in [m]: \cM_i \cap \cT \neq \emptyset\}$,
and note that $|\cB| \leq |\cT| \leq \overline{\cT}$ almost surely.
For $i \notin \cB$, we have $\wh\bU_i\wh\bU_i^\top = \bP_i^0$;
for $i \in \cB$, 
both $\wh\bU_i\wh\bU_i^\top$ and $\bP_i^0$ are still orthogonal projectors, so that
\[
- \bI_p \preceq 
\wh\bU_i \wh\bU_i^\top -\bP_i^0 
\preceq \bI_p ,
\qquad
\big\| \wh\bU_i \wh\bU_i^\top-\bP_i^0 \big\| \leq 1 ,
\]
which then implies
\begin{equation}
\label{eqn: diff_hatP_spectral}
\big\| \wh\bP - \wh\bP^0 \big\|
\leq 
\frac{1}{m} \sum_{i\in \cB}
\big\| \wh\bU_i \wh\bU_i^\top -\bP_i^0 \big\|
\leq 
\frac{|\cB|}{m}
\leq 
\frac{\overline{\cT}}{m} .
\end{equation}

From Lemma~\ref{lem: mean_block_fixed}, the matrix $\bP^\ast$ has leading $r$ eigenspace projector $\bPi_r^\ast$ and eigengap at least $(c_\ast/2) p^{- \max\{1 - \delta,0\}}$ for all sufficiently large $p$.
Then by the Davis--Kahan theorem \citep[Theorem~2]{yu2015useful}, we have
\begin{align}
\label{eq:intermediate}
    d(\wh\bA_{b,1:r},\bA_r^\ast)
    &\leq 
    \frac{4 \sqrt{2r}}{c_\ast} p^{\max\{1 - \delta,0\}}
    \big\| \wh\bP - \bP^\ast \big\| 
    \leq 
    \frac{4 \sqrt{2r}}{c_\ast} p^{\max\{1 - \delta,0\}}
    \left\{ \big\| \wh\bP^0 -\bP^\ast \big\|
    + \frac{\overline{\cT}}{m} \right\} ,
\end{align}
which, combined with 
\eqref{eqn: DC_ano_high_prob},
Lemma~\ref{lem: davis_kahan_intermediate_loading}, and the triangle inequality,
implies the desired high-probability bound in Theorem~\ref{thm: DC_ano_clean}.
The desired asymptotic rate of convergence is direct for fixed $b$,
and this concludes the proof.
\end{proof}


\begin{proof}[Proof of Theorem~\ref{thm: r_hat_consistency}]
If $r_{\max} = 2$, then $r=1$ and the maximization in~\eqref{eqn: r_hat_vfm} has only one candidate, and the conclusion trivially holds.
It hence suffices to consider $r_{\max} \geq 3$.
Throughout this proof, $r+1 \leq r_{\max} = b$ and $b = o(p)$.

Let $m_1 := \lfloor |\cI_1| / b \rfloor$ be the number of full blocks in $\cI_1$,
and by construction, $m_1 = m/2 + \cO(1)$.
As \eqref{eqn: thm_r_contam_condition} implies $\overline{\cT}/m \to 0$, we thus have
\[
m_1 - \overline{\cT} \asymp m ,
\qquad
m - \overline{\cT} \asymp m ,
\]
and in particular, $m_1 > \overline{\cT}$ for all sufficiently large $m,p$.
Assumption~\ref{ass: contamination} implies $|\cT \cap \cI_1| \leq \overline{\cT}$ almost surely.
Hence with the notation $r_m$ from Lemma~\ref{lem: stage2_variance}, we have
\[
\sqrt{\frac{\log(p)}{m_1 - \overline{\cT}}}
+ \frac{1}{p}
+ \frac{\overline{\cT}}{m_1}
=
\cO(r_m) .
\]
Separately, by the definition of $r_m$ and \eqref{eqn: thm_r_contam_condition}, we have
\begin{equation}
\label{eqn: star_rm_rate}
\sqrt{p} r_m^2
\leq
3 \left\{ \frac{\sqrt{p} \log(p)}{m - \overline{\cT}}
+ \frac{1}{p^{3/2}}
+ \sqrt{p}
\left( \frac{\overline{\cT}}{m} \right)^2
\right\}
\to 0 .
\end{equation}

For $j \in [r_{\max}]$, recall $\bb_j$ in Lemma~\ref{lem: stage2_variance}, and here, we can write $\bb_j= \bM^\top \wh\bA_{r_{\max}}^{(1)} \be_j$. 
Further define $B_p := \max_{r < j \leq r_{\max}} \|\bb_j\|^2$.
As \eqref{eqn: star_rm_rate} implies $r_m = o(1)$, we may apply Lemma~\ref{lem: stage2_variance},
so that there exist absolute constants $c_1, C_1 > 0$ such that, with probability approaching 1,
\begin{equation}
\label{eqn: b_j_norm_regimes}
c_1 p
\leq
\min_{j\in [r]} \|\bb_j\|^2
\leq
\max_{j\in [r]} \|\bb_j\|^2
\leq
C_1 p .
\end{equation}
The same lemma gives $\|\bb_j\|^2 \geq c_\epsilon$ for every $j \in [r_{\max}]$, 
and $c_\epsilon \leq B_p \leq C_\epsilon + \cO_P(p r_m^2)$.
Together with \eqref{eqn: star_rm_rate}, this implies
\begin{equation}
\label{eqn: star_Bp_rate}
\frac{B_p}{\sqrt{p}}
\leq
\frac{C_\epsilon}{\sqrt{p}}
+
\cO_P(\sqrt{p} \, r_m^2)
=
o_P(1) .
\end{equation}

Next, consider the conditional distribution of the clean ratios in $\cI_2$.
For clarity, define the sigma-algebra $\cF_1 := \sigma\big( \cT, \{\bx_t\}_{t\in \cI_1} \big)$.
The estimator $\wh\bA_{r_{\max}}^{(1)}$ and hence each $\bb_j$, is measurable with respect to $\cF_1$.
By Assumption~\ref{ass: sample_split_independence}
and independence of $\cT$ from the clean data, 
conditionally on $\cF_1$, the observations $\{\bx_t\}_{t \in\cI_2 \setminus \cT}$ remain independent with their original elliptical distribution.
Thus their conditional distribution admits the representation in~\eqref{eqn: xt_ellip_form}
(with independent copies of the original radius and uniform direction).
Then noting that $\P(\zeta > 0) = 1$ by Assumption~\ref{ass: elliptical},
and that $\bb_j \neq \zero$ by $\|\bb_j\|^2 \geq c_\epsilon$,
so that conditionally on $\cF_1$,
the clean denominators
in~\eqref{eqn: ratio_def} for $t \notin \cT$ are nonzero almost surely.
With \eqref{eqn: projected_coords}, the conditional distribution of a clean ratio (i.e., for $t\in \cI_2 \setminus \cT$) is hence
\[
\wh\rho_{t,j} =
\frac{\zeta_t^2 (\bb_j^\top \bv_t)^2}{\zeta_t^2 (\bb_{j+1}^\top \bv_t)^2}
=
\frac{(\bb_j^\top \bv_t)^2}{(\bb_{j+1}^\top \bv_t)^2} .
\]
Let $F_j$ denote the conditional distribution function of $\wh\rho_{t,j}$ for $t \in \cI_2 \setminus \cT$ given $\cF_1$.
By Lemma~\ref{lem: median_ratio_quadratic},
the population median of the conditional distribution $F_j$ is $m_j^\ast := \|\bb_j\|^2 /\|\bb_{j+1}\|^2$.
Note that on the event that \eqref{eqn: b_j_norm_regimes} holds,
we have
\begin{equation}
\label{eqn: pop_median_r}
m_r^\ast \geq
\frac{c_1 p}{B_p} ,
\end{equation}
and similarly,
\begin{equation*}
m_j^\ast \leq
\begin{cases}
C_1 /c_1 , & j < r , \\
B_p /c_\epsilon , & r < j \leq r_{\max}-1 ,
\end{cases}
\end{equation*}
which, together with $B_p \geq c_\epsilon$, implies that there exists an absolute constant $C_2 > 0$ such that $\max_{\substack{j\in [r_{\max}-1] \\ j \neq r}} m_j^\ast \leq C_2 B_p$ on the same event.
Combining this with \eqref{eqn: star_Bp_rate} and \eqref{eqn: pop_median_r}, we obtain
\begin{equation}
\label{eqn: pop_separation}
\frac{m_r^\ast}{\max_{\substack{j \in [r_{\max}-1] \\ j \neq r}} m_j^\ast}
\geq
\frac{c_1 p}{C_2 B_p^2}
\xrightarrow{P}
\infty .
\end{equation}

It remains to control the empirical medians, including the effect of anomalies in $\cI_2$.
To this end, define $\varepsilon_2 := |\cI_2 \cap \cT| /|\cI_2|$.
By Assumption~\ref{ass: contamination}, we have
\[
\varepsilon_2
\leq
\frac{\overline{\cT}}{|\cI_2|}
\leq
\frac{2 \overline{\cT}}{T}
\leq
\frac{2 \overline{\cT}}{m}
\to 0 ,
\]
almost surely,
and hence $|\cI_2 \setminus \cT| \geq |\cI_2| - \overline{\cT} \asymp T$.
As $\log(p) = o(m)$ by \eqref{eqn: thm_r_contam_condition}, 
we also have $\log(p) = o\big( |\cI_2| - \overline{\cT} \big)$.
Then for $j\in [r_{\max}-1]$, define the empirical distribution functions
\[
\wh{F}_j^{\mathrm{clean}}(s) :=
\frac{1}{|\cI_2 \setminus \cT|}
\sum_{t\in \cI_2 \setminus \cT} \mathbbm{1}_{\{\wh\rho_{t,j} \leq s\}} ,
\qquad
\wh{F}_j(s) :=
\frac{1}{|\cI_2|}
\sum_{t\in \cI_2} \mathbbm{1}_{\{\wh\rho_{t,j} \leq s\}} .
\]
Conditionally on $\cF_1$, the clean ratios $\wh{\rho}_{t, j}$ for $t \in \cI_2 \setminus \cT$
are i.i.d.\ across $t$,
and $m_j^\ast$ is fixed.
Hence by Lemma~\ref{lem: median_ratio_quadratic} with $\delta = 1/2$ and the notation $c_0 > 0$ therein, we have
\begin{equation}
\label{eqn: Fj_mj_two_bounds}
F_j\bigg( \frac{m_j^\ast}{2} \bigg)
\leq
\frac{1}{2} - \frac{c_0}{2} ,
\qquad
F_j\bigg( \frac{3 m_j^\ast}{2} \bigg)
\geq
\frac{1}{2} + \frac{c_0}{2} .
\end{equation}

Fix $j\in [r_{\max}-1]$ and $a \in \{1/2, 3/2\}$.
Then conditionally on $\cF_1$, the set $\cI_2 \setminus \cT$ and the threshold $a m_j^\ast$ are fixed,
together with the fact that $\{\wh\rho_{t,j}\}_{t\in \cI_2 \setminus \cT}$ are i.i.d.\ under the conditional distribution,
it holds that $\big\{ \mathbbm{1}_{\{\wh\rho_{t,j} \leq a m_j^\ast\}} \big\}_{t\in \cI_2 \setminus \cT}$ are conditionally independent,
and satisfy $\E\big[ \mathbbm{1}_{\{\wh\rho_{t,j} \leq a m_j^\ast\}} \mid \cF_1 \big] = F_j(a m_j^\ast)$.
Hence by the Hoeffding's inequality \citep[Theorem~2.2.6]{vershynin2018}, for every $u > 0$, we have
\begin{align*}
    &\;\quad
    \P\Big( \big| \wh{F}_j^{\mathrm{clean}}(a m_j^\ast)
    - F_j(a m_j^\ast) \big|
    > u \mid \cF_1 \Big) \\
    &=
    \P\Bigg( \bigg| \frac{1}{|\cI_2 \setminus \cT|}
    \sum_{t\in \cI_2 \setminus \cT}
    \big\{ \mathbbm{1}_{\{\wh\rho_{t,j} \leq a m_j^\ast\}}
    - F_j(a m_j^\ast)
    \big\} \bigg|
    > u \mid \cF_1 \Bigg) \\
    &=
    \P\Bigg( \bigg|
    \sum_{t\in \cI_2 \setminus \cT}
    \big\{ \mathbbm{1}_{\{\wh\rho_{t,j} \leq a m_j^\ast\}}
    - \E\big[ \mathbbm{1}_{\{\wh\rho_{t,j} \leq a m_j^\ast\}} \mid \cF_1 \big]
    \big\} \bigg|
    > u \cdot |\cI_2 \setminus \cT| \mid \cF_1 \Bigg) 
    \leq
    2 \exp\big( -2 |\cI_2 \setminus \cT| \cdot u^2 \big) ,
\end{align*}
which, by taking $u = c_0/8$
and a union bound argument, gives
\begin{align}
    &\;\quad
    \P\left( \max_{j\in [r_{\max}-1]}
    \max_{a \in \{1/2, 3/2\}}
    \big| \wh{F}_j^{\mathrm{clean}}(a m_j^\ast) - F_j(a m_j^\ast) \big|
    > \frac{c_0}{8}
    \mid \cF_1 \right) \notag \\
    &\leq
    \sum_{j=1}^{r_{\max}-1}
    \sum_{a\in \{1/2, 3/2\}}
    \P\left( \big| \wh{F}_j^{\mathrm{clean}}(a m_j^\ast) - F_j(a m_j^\ast) \big|
    > \frac{c_0}{8}
    \mid \cF_1 \right) \notag \\
    &\leq
    4 (r_{\max}-1)
    \exp\left( - \frac{|\cI_2 \setminus \cT| \cdot c_0^2}{32} \right) 
    \leq
    4p \exp\left( - \frac{(|\cI_2| - \overline{\cT}) \cdot c_0^2}{32} \right)
    = o(1) ,
    \label{eqn: diff_ecdf_cdf_clean}
\end{align}
uniformly over the conditioning,
where the last inequality used $|\cI_2 \setminus \cT| \geq |\cI_2| - \overline{\cT}$ 
and $r_{\max} = b = o(p)$,
and the last equality used $\log(p) = o(|\cI_2| - \overline{\cT})$ by $\log(p) = o(m)$, $m\leq T$, and $T\asymp |\cI_2| - \overline{\cT}$.
Note that by taking expectations of the conditional probability bound, \eqref{eqn: diff_ecdf_cdf_clean} also holds unconditionally.

For any finite $s$, we have the deterministic inequality
\[
\big| \wh{F}_j(s) - F_j(s) \big|
\leq
(1-\varepsilon_2) \cdot \big| \wh{F}_j^{\mathrm{clean}}(s) - F_j(s) \big|
+ \varepsilon_2 
\leq
\big| \wh{F}_j^{\mathrm{clean}}(s) - F_j(s) \big| 
+ \varepsilon_2 ,
\]
which, noting that $\varepsilon_2 \to 0$ so that eventually $\varepsilon_2 \leq c_0/8$,
together with \eqref{eqn: Fj_mj_two_bounds} 
and the unconditional version of \eqref{eqn: diff_ecdf_cdf_clean} justified above,
it holds with probability approaching 1 that
\[
\wh{F}_j\bigg( \frac{m_j^\ast}{2} \bigg)
\leq
F_j\bigg( \frac{m_j^\ast}{2} \bigg)
+ \frac{c_0}{4}
\leq
\frac{1}{2} - \frac{c_0}{4} ,
\qquad
\wh{F}_j\bigg( \frac{3 m_j^\ast}{2} \bigg)
\geq
F_j\bigg( \frac{3 m_j^\ast}{2} \bigg)
- \frac{c_0}{4}
\geq
\frac{1}{2} + \frac{c_0}{4} ,
\]
so that with the fact that $1/2 - c_0/4 < \lceil |\cI_2|/2 \rceil /|\cI_2| \leq 1/2 + c_0/4$ for sufficiently large $|\cI_2|$, we obtain
\begin{equation}
\label{eqn: ecdf_upper_lower}
\wh{F}_j\bigg( \frac{m_j^\ast}{2} \bigg)
<
\frac{\lceil |\cI_2|/2 \rceil}{|\cI_2|}
\leq
\wh{F}_j\bigg( \frac{3 m_j^\ast}{2} \bigg) ,
\end{equation}
uniformly for all $j\in [r_{\max}-1]$.

Recall that $\wh\rho_j$ in~\eqref{eqn: median_ratio} is the $\lceil |\cI_2|/2 \rceil$th smallest ratio,
and also recall the definition that $\wh{F}_j(s) = |\cI_2|^{-1} \sum_{t\in \cI_2} \mathbbm{1}_{\{\wh\rho_{t,j} \leq s\}}$.
Together with \eqref{eqn: ecdf_upper_lower}, we finally conclude
\begin{equation*}
\frac{m_j^\ast}{2}
<
\wh\rho_j
\leq
\frac{3 m_j^\ast}{2} ,
\qquad
j \in [r_{\max}-1] ,
\end{equation*}
simultaneously with probability approaching 1,
which, combined with \eqref{eqn: pop_separation},
implies that on an event with probability approaching 1,
\[
\frac{\wh\rho_r}{\max_{\substack{j \in [r_{\max}-1] \\ j \neq r}} \wh\rho_j}
\geq
\frac{1}{3} \cdot
\frac{m_r^\ast}{\max_{\substack{j \in [r_{\max}-1] \\ j \neq r}} m_j^\ast} 
\xrightarrow{P}
\infty .
\]
By the definition in~\eqref{eqn: r_hat_vfm}, this concludes $\P(\wh{r} = r) \to 1$ and thus completes this proof.
\end{proof}

\clearpage

\subsection{Auxiliary results with proofs}
\label{subsec: lemmas}

\subsubsection{Lemmas for DC}

\begin{lemma}[Uniqueness of $\wh\bU_i \wh\bU_i^\top$]
\label{lem: block_projector_unique}
Let Assumptions~\ref{ass: elliptical}, \ref{ass: noise},
and \ref{ass: contamination} hold,
and suppose $1 \leq b \leq p$.
Conditional on $\cT$,
for any clean block~$i\in[m]$ with $\cM_i \cap \cT = \emptyset$ (i.e., there are no anomalies in the $i$th interval),
the matrix $\wh\bS_i$ in~\eqref{eqn: block_cov_vfm} has rank $b$ almost surely.
Hence its leading $b$ eigenspace has a uniquely defined orthogonal projector $\wh\bU_i \wh\bU_i^\top$,
also almost surely.
\end{lemma}

\begin{proof}[Proof of Lemma~\ref{lem: block_projector_unique}]
Condition on $\cT$ and fix any clean block~$i\in[m]$.
For simplicity, let $\cM_i$ be explicitly $\cM_i = \{t_1, \dots, t_b\}$,
and define $\bX_i := [\bx_{t_1}, \dots, \bx_{t_b}] \in \R^{p \times b}$,
$\bY_i := [\bM \bu_{t_1}, \dots, \bM \bu_{t_b}]$,
and $\bD_i := \diag\big( \zeta_{t_1} /\|\bu_{t_1}\|, \dots, \zeta_{t_b} /\|\bu_{t_b}\| \big)$.
Then $\bX_i = \bY_i \bD_i$ by \eqref{eqn: xt_ellip_form}, 
and $\bD_i$ is invertible almost surely by Assumption~\ref{ass: elliptical}.

As we also have $\bSigma_x \succeq c_\epsilon \bI_p$ by Assumption~\ref{ass: noise} and \eqref{eqn: cov},
we may conclude that the columns of $\bY_i$ are independent Gaussian vectors with covariance matrix $\bSigma_x$,
and hence these columns are linearly independent almost surely as $b \leq p$.
Thus, it holds almost surely that the rank of $\bY_i$ is $b$, 
and so is the rank of $\bX_i$ and hence that of $\wh\bS_i$ by $\wh\bS_i = b^{-1} \bX_i \bX_i^\top$.
Namely, it holds almost surely that $\wh\bS_i$ has $b$ positive eigenvalues and $p-b$ zero eigenvalues.
Hence its leading $b$ eigenspace, 
i.e., $\wh\bU_i \wh\bU_i^\top$, is uniquely defined almost surely,
which concludes this proof.
\end{proof}

\begin{lemma}[Concentration of the aggregated projector]
\label{lem: subgauss_aggregate_P}
Let Assumptions~\ref{ass: elliptical} and \ref{ass: noise} hold,
and suppose there are no anomalies, i.e., $\cT=\emptyset$ in~\eqref{eqn: contam_index}.
Suppose $1 \leq b \leq p$.
Then the aggregated projector in~\eqref{eqn: aggregate_P} satisfies that,
for every $x \in (0,1)$,
\begin{align*}
    \P\left( \big\|\wh\bP -\E[\wh\bP] \big\| >
    4 \sqrt{\frac{2 \log(2p/x)}{m}} 
    \right) \leq  x 
    &, \\
    \P\left( \big\| \wh\bP - \E[\wh\bP] \big\| >
    \sqrt{\frac{b}{m}}
    + \sqrt{\frac{\log(1/x)}{2m}}
    \right) \leq  x 
    &.
\end{align*}
\end{lemma}

\begin{proof}[Proof of Lemma~\ref{lem: subgauss_aggregate_P}]
By Lemma~\ref{lem: block_projector_unique},
each $\wh\bU_i \wh\bU_i^\top$ is uniquely defined with probability~1.
Hence using notations in~\eqref{eqn: block_cov_vfm}, for each block~$i\in[m]$,
define $\wh{\bPi}_i := \wh{\bU}_i \wh{\bU}_i^\top -\E[\wh{\bU}_i \wh{\bU}_i^\top]$.
By Assumption~\ref{ass: elliptical}, $\{\wh{\bPi}_i\}_{i=1}^m$ are i.i.d.\ symmetric random matrices in $\R^{p \times p}$ with $\E[\wh{\bPi}_i] = \zero$.
Also, note that by construction and by the triangle inequality, for every $i\in[m]$ we have
\begin{align*}
    \|\wh{\bPi}_i\| \leq
    \|\wh{\bU}_i \wh{\bU}_i^\top\| + \|\E[\wh{\bU}_i \wh{\bU}_i^\top]\|
    \leq 2 .
\end{align*}
Hence the matrix Azuma's inequality \citep[e.g.][Theorem~7.1]{Tropp2012} is applicable, which gives that, every $u > 0$,
\begin{equation}\label{eqn: mat_Azuma}
\P\left( \left\| \frac{1}{m} \sum_{i = 1}^{m} \wh{\bPi}_i \right\| \geq u \right) 
\leq 
2p \exp\left( -\frac{m u^{2}}{32} \right) .
\end{equation}
Setting $u = 4 \sqrt{2 \log(2p /x) /m}$ gives $m u^2 /32 = \log(2p /x)$, so that the right-hand side of \eqref{eqn: mat_Azuma} is $x$.
The proof of the first bound is then complete by noting that by \eqref{eqn: aggregate_P},
\[
\left\|\wh\bP -\E[\wh\bP] \right\| 
=
\left\| \frac{1}{m} \sum_{i = 1}^{m} \big\{ \wh{\bU}_i \wh{\bU}_i^\top - \E[\wh{\bU}_i \wh{\bU}_i^\top] \big\} \right\| 
=
\left\| \frac{1}{m} \sum_{i = 1}^{m} \wh{\bPi}_i \right\| .
\]

For the second bound, similarly by Lemma~\ref{lem: block_projector_unique}, the matrices $\wh\bU_i\wh\bU_i^\top$ are independent rank-$b$ orthogonal projectors almost surely,
with $\|\wh\bU_i\wh\bU_i^\top\|_F^2 = \tr(\wh\bU_i\wh\bU_i^\top) = b$.
Hence by Jensen's inequality, we have
\begin{align}
    &\;\quad
    \E\left[ \big\| \wh\bP - \E[\wh\bP] \big\| \right]
    \leq
    \E\left[ \big\| \wh\bP - \E[\wh\bP] \big\|_F \right] 
    =
    \E\left[ \left\| \frac{1}{m} \sum_{i=1}^m \big( \wh\bU_i\wh\bU_i^\top - \E[\wh\bU_i\wh\bU_i^\top] \big) \right\|_F \right] \notag \\
    &\leq
    \left\{ \E\left[ \left\| \frac{1}{m} \sum_{i=1}^m \big( \wh\bU_i\wh\bU_i^\top - \E[\wh\bU_i\wh\bU_i^\top] \big) \right\|_F^2 \right] \right\}^{1/2} 
    =
    \left\{ \frac{1}{m^2} \sum_{i=1}^m \E\big[ \big\| \wh\bU_i\wh\bU_i^\top - \E[\wh\bU_i\wh\bU_i^\top] \big\|_F^2 \big] \right\}^{1/2} \notag \\
    &=
    \left\{ \frac{1}{m} \big\{ \E\big[ \|\wh\bU_i\wh\bU_i^\top\|_F^2 \big] - \big\| \E[\wh\bU_i\wh\bU_i^\top] \big\|_F^2 \big\} \right\}^{1/2}
    \leq
    \sqrt{\frac{b}{m}} .
    \label{eqn: hatP_demean_mean}
\end{align}

Next, consider replacing $\wh\bU_i\wh\bU_i^\top$ by another rank-$b$ orthogonal projector $\wh\bU_i' (\wh\bU_i')^\top$,
and note that
\[
-\bI_p
\preceq
\wh\bU_i\wh\bU_i^\top
-
\wh\bU_i' (\wh\bU_i')^\top
\preceq
\bI_p ,
\]
which implies $\big\| \wh\bU_i\wh\bU_i^\top - \wh\bU_i' (\wh\bU_i')^\top \big\| \leq 1$.
By the reverse triangle inequality, replacing this block changes the value of $\big\| \wh\bP - \E[\wh\bP] \big\|$ by at most $1/m$.
Together with the fact that $\wh\bU_i\wh\bU_i^\top$ are independent with each other,
we may apply McDiarmid's inequality \citep[][Theorem~2.9.1]{vershynin2018},
with bounded differences $1/m$ for each of the $m$ blocks,
so that for every $u > 0$,
\[
\P\left(
\big\| \wh\bP - \E[\wh\bP] \big\|
- \E\left[ \big\| \wh\bP - \E[\wh\bP] \big\| \right]
> u \right)
\leq
\exp\left( - \frac{2u^2}{\sum_{i=1}^m (1/m)^2} \right)
=
\exp(-2 m u^2) ,
\]
which, combined with \eqref{eqn: hatP_demean_mean}, gives
\[
\P\left(
\big\| \wh\bP - \E[\wh\bP] \big\|
> \sqrt{\frac{b}{m}} + u
\right)
\leq
\exp(-2 m u^2) .
\]
Setting $u = \sqrt{\log(1/x)/(2m)}$ proves the second bound
and thus completes the proof.
\end{proof}

\begin{lemma}
\label{lem: davis_kahan_intermediate_loading}
Let Assumptions~\ref{ass: elliptical}, \ref{ass: loading}, and \ref{ass: noise} hold.
Recall that $\ddot\bA \in\R^{p\times r}$ consists of an orthonormal basis of $\spa(\bA)$,
and $\bA_r^\ast = \mathrm{SVD}_r(\bSigma_x)$.
Then for all sufficiently large $p$,
we have
\[
d(\bA_r^\ast,\ddot\bA) 
\leq
\frac{2 \sqrt{2r} C_\epsilon}{c_\lambda p^\delta} .
\]
\end{lemma}

\begin{proof}[Proof of Lemma~\ref{lem: davis_kahan_intermediate_loading}]
By \eqref{eqn: cov} and Assumption~\ref{ass: noise}, 
we have $\big\| \bSigma_x -  \bA\bA^\top \big\| = \big\| \bSigma_\epsilon \big\| \leq C_\epsilon$.
Also, the gap between the $r$th and $(r + 1)$th eigenvalues of $\bA\bA^\top$ is $\lambda_r(\bA^\top\bA) \geq c_\lambda p^\delta$ by Assumption~\ref{ass: loading}.
Hence by the Davis--Kahan theorem \citep[Theorem~2]{yu2015useful}, we have
\[
d(\bA_r^\ast,\ddot\bA)
\leq
\frac{2 \sqrt{2r} \big\| \bSigma_x -\bA\bA^\top \big\|}{\lambda_r(\bA^\top\bA)}
\leq
\frac{2 \sqrt{2r} C_\epsilon}{c_\lambda p^{\delta}} ,
\]
as desired, which completes this proof.
\end{proof}


\begin{lemma}
\label{lem: mean_block_fixed}
Let Assumptions~\ref{ass: elliptical}, \ref{ass: loading}, \ref{ass: noise},
and \ref{ass: contamination} hold.
Suppose $r \leq b$ and $b = o(p^{\min\{\delta,1\}})$.
Recall $\bV_x$, $\bLambda_x$, $\bA_r^\ast$, and $\bPi_r^\ast$
defined at the beginning of this supplement.
Then 
for every clean block~$i\in [m]$ (i.e., there are no anomalies in the $i$th block),
there exist $\gamma_1, \dots, \gamma_p$ such that $$
\E[\wh\bU_i\wh\bU_i^\top \mid\cT] = \bV_x \diag(\gamma_1, \dots, \gamma_p) \bV_x^\top.
$$

In particular, with $d := b - r + 1$, for all sufficiently large $p$, it holds that
\begin{equation}
\label{eqn: gamma_j_upper_lower}
\begin{aligned}
\gamma_j &\geq
\frac{d\lambda_j(\bSigma_x)}{\sum_{k=r+1}^p \lambda_k(\bSigma_x) + (d+2) \lambda_j(\bSigma_x)} ,
&& 
j\in[r] 
,\\
\gamma_j &\leq
\frac{b \lambda_j(\bSigma_x)}{c_\epsilon (p-r)} ,
&&
r< j \leq p .
\end{aligned}
\end{equation}
Also, with $c_\ast := c_\lambda /(C_\epsilon + 3 c_\lambda) >0$, it holds that
\[
\min_{j\in[r]} \gamma_j
\geq 
c_\ast p^{- \max\{1 - \delta,0\}},
\qquad
\max_{r < j \leq p} \gamma_j
\leq
\frac{b C_\epsilon}{c_\epsilon (p-r)} ,
\]
and the leading $r$ eigenspace of $\E[\wh\bU_i\wh\bU_i^\top \mid\cT]$ has orthogonal projector $\bPi_r^\ast$, with the corresponding eigengap
$\lambda_r\big( \E[\wh\bU_i\wh\bU_i^\top \mid\cT] \big) - \lambda_{r+1}\big( \E[\wh\bU_i\wh\bU_i^\top \mid\cT] \big) \geq (c_\ast/2) p^{- \max\{1 - \delta,0\}}$ for all sufficiently large $p$.
\end{lemma}

\begin{proof}[Proof of Lemma~\ref{lem: mean_block_fixed}]
Condition on $\cT$ and fix a clean block $i\in [m]$.
For $t\in \cM_i$, define $\bg_t := \bLambda_x^{-1/2} \bV_x^\top \bM \bu_t$.
Then by \eqref{eqn: cov}, we can write
\[
\E[\bg_t \bg_t^\top] =
\bLambda_x^{-1/2}
\bV_x^\top \bM \bM^\top \bV_x
\bLambda_x^{-1/2}
= \bI_p .
\]
Hence the vectors $\bg_t$ are independent standard Gaussian random vectors in $\R^p$.
Recall the notations $\cM_i = \{t_1, \dots, t_b\}$, $\bX_i$, and $\bD_i$ in the proof of Lemma~\ref{lem: block_projector_unique}.
Also define $\bG_i := [\bg_{t_1}, \dots, \bg_{t_b}] \in \R^{p \times b}$,
and denote its $j$th row as $\bg_{i,j}^\top$, where $\bg_{i,j} \in \R^b$.
Then by construction, $\bV_x^\top \bX_i = \bLambda_x^{1/2} \bG_i \bD_i$ 
and since $\bD_i$ is invertible almost surely under Assumption~\ref{ass: elliptical}, by Lemma~\ref{lem: block_projector_unique} we have
\begin{equation}
\label{eqn: VUUV_equal}
\bV_x^\top \wh\bU_i \wh\bU_i^\top \bV_x
=
\bLambda_x^{1/2} \bG_i
(\bG_i^\top \bLambda_x \bG_i)^{-1}
\bG_i^\top \bLambda_x^{1/2} ,
\end{equation}
where the inverse exists almost surely as $\bG_i$ has full column rank and $\bLambda_x \succeq c_\epsilon \bI_p$.

For any $j \in [p]$, let $\bD_{e,j} := \bI_p - 2 \be_j \be_j^\top$,
and the matrices $\bD_{e,j} \bG_i$ and $\bG_i$ have the same distribution.
As $\bD_{e,j} \bLambda_x \bD_{e,j} = \bLambda_x$,
following the proof of Theorem~2 of \cite{Fanetal2019},
by \eqref{eqn: VUUV_equal} we have
\[
\bV_x^\top
\E[\wh\bU_i \wh\bU_i^\top \mid \cT]
\bV_x
=
\bD_{e,j}
\bV_x^\top
\E[\wh\bU_i \wh\bU_i^\top \mid \cT]
\bV_x
\bD_{e,j} ,
\]
which holds for any $j \in [p]$, 
and hence $\bV_x^\top \E[\wh\bU_i \wh\bU_i^\top \mid \cT] \bV_x$ must be diagonal with $\E[\wh\bU_i \wh\bU_i^\top \mid \cT] = \bV_x \diag(\gamma_1, \dots, \gamma_p) \bV_x^\top$,
where $\gamma_j = \lambda_j(\bSigma_x) \E\big[ \bg_{i,j}^\top (\bG_i^\top \bLambda_x \bG_i)^{-1} \bg_{i,j} \big]$.

Fix $j\in [r]$. Let $\bR_j \in \R^{b\times b}$ be the orthogonal projection matrix onto the orthogonal complement of $\spa\{\bg_{i,k}: k\in [r], k \neq j\}$,
where we set $\bR_j = \bI_b$ in the case of $r = 1$.
Then note that $\bR_j$ has rank $d = b -r +1$, and is independent of $\bg_{i,j}$.
Define $\bv:= \bR_j\bg_{i,j}$, 
$s:= \bg_{i,j}^\top \bR_j \bg_{i,j}$,
and $\bH_0 := \sum_{k=r+1}^p \lambda_k(\bSigma_x) \bg_{i,k} \bg_{i,k}^\top$.
Then conditional on $\bR_j$, the variable $s \sim \chi_d^2$, so that
\begin{equation}
\label{eqn: g_moments}
\E[s]=d,
\qquad
\E[s^2]=d(d+2),
\qquad
\|\bv\|^2=s.
\end{equation}
We also have
\begin{align}
    \E[\bv^\top \bH_0 \bv] &=
    \sum_{k=r+1}^p \lambda_k(\bSigma_x) 
    \cdot \E\left[ \E\big\{ (\bg_{i,k}^\top \bv)^2 \mid \bv \big\} \right] 
    =
    \sum_{k=r+1}^p \lambda_k(\bSigma_x) \E\big[ \|\bv\|^2 \big]
    =
    d \sum_{k=r+1}^p \lambda_k(\bSigma_x) .
    \label{eqn: E_vH0v}
\end{align}

As $\bG_i^\top \bLambda_x \bG_i$ is positive definite and $\bv$ is nonzero, by the Cauchy--Schwarz inequality, we have
\begin{align*}
    s^2 &= \big( \bg_{i,j}^\top \bv \big)^2
    =
    \Big\{ \big( (\bG_i^\top \bLambda_x \bG_i)^{-1/2} \bg_{i,j} \big)^\top 
    \big( (\bG_i^\top \bLambda_x \bG_i)^{1/2} \bv \big) 
    \Big\}^2 \\
    &\leq
    \big( \bg_{i,j}^\top (\bG_i^\top \bLambda_x \bG_i)^{-1} \bg_{i,j} \big) 
    \big( \bv^\top \bG_i^\top \bLambda_x \bG_i \bv \big) ,
\end{align*}
which, together with the fact that $\bv^\top \bG_i^\top \bLambda_x \bG_i \bv = \lambda_j(\bSigma_x) s^2 + 0 + \bv^\top \bH_0 \bv$,
implies
\[
\lambda_j(\bSigma_x) \cdot \bg_{i,j}^\top (\bG_i^\top \bLambda_x \bG_i)^{-1} \bg_{i,j}
\geq
\frac{\lambda_j(\bSigma_x) \cdot s^2}{\lambda_j(\bSigma_x) s^2 + \bv^\top \bH_0 \bv} .
\]
Combining this with \eqref{eqn: g_moments} and \eqref{eqn: E_vH0v}, we obtain
\begin{align*}
    \gamma_j
    &=
    \E\left[ \lambda_j(\bSigma_x) \cdot \bg_{i,j}^\top (\bG_i^\top \bLambda_x \bG_i)^{-1} \bg_{i,j} \right]
    \geq
    \E\left[ \frac{s^2}{s^2 + \bv^\top \bH_0 \bv /\lambda_j(\bSigma_x)} \right] \\
    &\geq
    \frac{\{\E[s]\}^2}{\E\big[ s^2 + \bv^\top \bH_0 \bv /\lambda_j(\bSigma_x) \big]}
    =
    \frac{\{\E[s]\}^2}{\E[s^2] + \E[\bv^\top \bH_0 \bv] /\lambda_j(\bSigma_x)} \\
    &=
    \frac{d^2}{d (d+2) + d \sum_{k=r+1}^p \lambda_k(\bSigma_x) /\lambda_j(\bSigma_x)} 
    =
    \frac{d \lambda_j(\bSigma_x)}{\sum_{k=r+1}^p \lambda_k(\bSigma_x) + (d+2) \lambda_j(\bSigma_x)} ,
\end{align*}
which is the desired lower bound in~\eqref{eqn: gamma_j_upper_lower},
where the second inequality used the Cauchy--Schwarz inequality on $s = (s/\sqrt{s^2 + \bv^\top \bH_0 \bv /\lambda_j(\bSigma_x)}) \sqrt{s^2 + \bv^\top \bH_0 \bv /\lambda_j(\bSigma_x)}$.

For the upper bound in~\eqref{eqn: gamma_j_upper_lower}, let $\bG_{i,0} \in \R^{(p-r) \times b}$ denote the sub-matrix containing rows $r+1, \dots, p$ of $\bG_i$.
As $b = o(p^{\min\{\delta,1\}})$, we have $b \leq p - r$ for all sufficiently large $p$.
Hence $\bG_{i,0}^\top \bG_{i,0}$ is invertible almost surely, and
\[
\bG_i^\top \bLambda_x \bG_i 
\succeq
c_\epsilon \bG_{i,0}^\top \bG_{i,0}
, 
\qquad
(\bG_i^\top \bLambda_x \bG_i)^{-1} 
\preceq 
c_\epsilon^{-1} (\bG_{i,0}^\top \bG_{i,0})^{-1} .
\]
Note that the matrix $\bG_{i,0}(\bG_{i,0}^\top \bG_{i,0})^{-1} \bG_{i,0}^\top$ is a rank-$b$ projection matrix with exchangeable diagonal entries
and trace exactly $b$, 
so that each of its diagonal entries has expectation $b/(p - r)$.
Hence for $j>r$, we have
\[
\gamma_j
\leq
\frac{\lambda_j(\bSigma_x)}{c_\epsilon}
\E\big[ \bg_{i,j}^\top (\bG_{i,0}^\top \bG_{i,0})^{-1} \bg_{i,j} \big]
=
\frac{b \lambda_j(\bSigma_x)}{c_\epsilon (p-r)} ,
\]
which is the desired result in~\eqref{eqn: gamma_j_upper_lower}.

Finally, for $j \leq r$, since $d \geq 1$ and $\sum_{k=r+1}^p \lambda_k(\bSigma_x) \geq 0$, we have
\begin{align*}
    \gamma_j
    &\geq
    \frac{\lambda_j(\bSigma_x)}{\sum_{k=r+1}^p \lambda_k(\bSigma_x) + 3 \lambda_j(\bSigma_x)}
    \geq
    \frac{c_\lambda p^\delta}{C_\epsilon p + 3 c_\lambda p^\delta} 
    \geq
    \frac{c_\lambda}{C_\epsilon + 3 c_\lambda}
    \min\{1, p^{\delta-1}\}
    =
    c_\ast p^{- \max\{1 - \delta,0\}} ,
\end{align*}
whereas for $j>r$, the upper bound in~\eqref{eqn: gamma_j_upper_lower} and $\lambda_j(\bSigma_x) \leq C_\epsilon$ apply.
Moreover, we have
\[
\frac{b /(p-r)}{p^{- \max\{1 - \delta,0\}}}
\asymp
\frac{b}{p^{\min\{\delta,1\}}}
\to
0 ,
\]
so that $\min_{j\in[r]} \gamma_j - \max_{r<j\leq p} \gamma_j \geq (c_\ast/2) p^{- \max\{1 - \delta,0\}}$ for all sufficiently large $p$,
which thus concludes the proof.
\end{proof}


\subsubsection{Lemmas for STAR}

\begin{lemma}
\label{lem: stage2_variance}
Let Assumptions~\ref{ass: elliptical},
\ref{ass: loading}, \ref{ass: noise},
\ref{ass: contamination},
and \eqref{ass: loading_star} hold.
Suppose $r< r_{\max} = b$ and $b = o(p)$.
Recalling $\wh\bA_b$ from \eqref{eqn: aggregate_P}, 
denote by $\wh\bA_{r_{\max}} := \wh\bA_b$ for the choice of $b = r_{\max}$,
and let $\wh\bA_{r_{\max},1:r}$ denote its first $r$ columns.
For each $j\in [r_{\max}]$, let $\be_j$ denote the $j$th standard basis vector in $\R^{r_{\max}}$, and define
\[
\bb_j := \bM^{\top} \wh\bA_{r_{\max}} \be_j \in \R^{r+p} ,
\qquad
r_m := \sqrt{\frac{\log(p)}{m - \overline{\cT}}} + \frac{1}{p} + \frac{\overline{\cT}}{m} .
\]
Also suppose $\min(m,p) \to \infty$, $m > \overline{\cT}$, and $r_m = o(1)$.
Then there exist absolute constants $c_1, C_1 > 0$ such that
\[
\P\left(
c_1 p
\leq
\min_{j\in [r]} \|\bb_j\|^2
\leq
\max_{j\in [r]} \|\bb_j\|^2
\leq
C_1 p
\right)
\to 1 .
\]
Moreover, with the notations in Assumption~\ref{ass: noise},
$\|\bb_j\|^2 \geq c_\epsilon$
for every $j \in [r_{\max}]$, and if $r_{\max} > r$,
\[
\max_{r < j \leq r_{\max}} \|\bb_j\|^2
\leq
C_\epsilon + \cO_P(p r_m^2) .
\]
\end{lemma}

\begin{proof}[Proof of Lemma~\ref{lem: stage2_variance}]
Under \eqref{ass: loading_star}, Assumption~\ref{ass: loading} holds with $\delta=1$.
Using \eqref{eq:intermediate} with $b=r_{\max}$ and Lemma~\ref{lem: subgauss_aggregate_P}, we have
\[
\left\| \wh\bA_{r_{\max},1:r} \wh\bA_{r_{\max},1:r}^\top
- \bPi_r^\ast \right\|_F
=
\cO_P\left( \sqrt{\frac{\log(p)}{m}}
+ \frac{\overline{\cT}}{m} \right) 
=
\cO_P(r_m) ,
\]
where the last equality used $m > \overline{\cT}$.
Hence
\begin{equation}
\label{eqn: dc_top_r_rate}
\left\| \wh\bA_{r_{\max},1:r} \wh\bA_{r_{\max},1:r}^\top
- \bPi_r^\ast \right\|_F
=
\cO_P(r_m) .
\end{equation}
Note that by \eqref{eqn: cov} and the definition of $\bb_j$, we can write $\|\bb_j\|^2 = \big( \wh\bA_{r_{\max}} \be_j \big)^\top \bSigma_x \big( \wh\bA_{r_{\max}} \be_j \big)$.
Since $\|\wh\bA_{r_{\max}} \be_j\| = 1$ and $\bSigma_x \succeq c_\epsilon \bI_p$ by Assumption~\ref{ass: noise}, 
we arrive at $\|\bb_j\|^2 \geq c_\epsilon$ for every $j \in [r_{\max}]$.

Consider $j \leq r$.
Since $\wh\bA_{r_{\max}} \be_j$ is in the column space of $\wh\bA_{r_{\max},1:r}$, we have
\begin{align*}
    \big\| (\bI_p - \bPi_r^\ast) \wh\bA_{r_{\max}} \be_j \big\| 
    &=
    \big\| \big( \wh\bA_{r_{\max},1:r} \wh\bA_{r_{\max},1:r}^\top - \bPi_r^\ast \big) \wh\bA_{r_{\max}} \be_j \big\| \\
    &\leq
    \big\| \wh\bA_{r_{\max},1:r} \wh\bA_{r_{\max},1:r}^\top - \bPi_r^\ast \big\|_F
    =
    \cO_P(r_m) ,
\end{align*}
uniformly over $j\in [r]$,
where the last equality used \eqref{eqn: dc_top_r_rate}.
Hence with the fact that $\bPi_r^\ast \wh\bA_{r_{\max}} \be_j$ is orthogonal to $(\bI_p - \bPi_r^\ast) \wh\bA_{r_{\max}} \be_j$, we obtain
\begin{equation}
\label{eqn: j_leq_r_PiAe}
\min_{j\in [r]} 
\big\| \bPi_r^\ast \wh\bA_{r_{\max}} \be_j \big\|^2
=
\min_{j\in [r]} \Big\{
\big\| \wh\bA_{r_{\max}} \be_j \big\|^2
- \big\| (\bI_p - \bPi_r^\ast) \wh\bA_{r_{\max}} \be_j \big\|^2
\Big\}
\geq
1 - \cO_P(r_m^2) .
\end{equation}
Note that it is direct that $\max_{j\in [r]}  \big\| \bPi_r^\ast \wh\bA_{r_{\max}} \be_j \big\|^2 \leq 1$.

By \eqref{ass: loading_star}
and Assumption~\ref{ass: noise}, we have $\lambda_r(\bSigma_x) \geq c_\lambda p$ and $\lambda_1(\bSigma_x) \leq C_\lambda p + C_\epsilon$.
Then by \eqref{eqn: j_leq_r_PiAe}, we have
\begin{align*}
    \min_{j\in [r]} \|\bb_j\|^2
    &\geq
    \lambda_r(\bSigma_x)
    \min_{j\in [r]} \big\| \bPi_r^\ast \wh\bA_{r_{\max}} \be_j \big\|^2 
    \geq
    c_\lambda p
    \big\{ 1 - \cO_P(r_m^2) \big\} ,
\end{align*}
and also $\max_{j\in [r]} \|\bb_j\|^2 \leq \lambda_1(\bSigma_x) \leq C_\lambda p + C_\epsilon$.
Together with $r_m = o(1)$,
these conclude the proof for the results for $j \leq r$.

It remains to show the upper bound for $j > r$ (the lower bound has been established at the beginning of this proof).
First, we decompose
\[
\wh\bA_{r_{\max}}\be_j = 
\bPi_r^\ast \wh\bA_{r_{\max}}\be_j + (\bI_p - \bPi_r^\ast) \wh\bA_{r_{\max}}\be_j 
=: \bs_j + \bn_j ,
\]
where, note that immediately, $\bs_j$ is orthogonal to $\bn_j$,
and hence $1= \| \wh\bA_{r_{\max}}\be_j \|^2 = \|\bs_j\|^2 + \|\bn_j\|^2$.
For $\bs_j$, noting that the vector $\wh\bA_{r_{\max}} \be_j$ is orthogonal to $\wh\bA_{r_{\max},1:r}$ by construction, we have
\begin{equation}
\label{eqn: s_j_rate}
\|\bs_j\| = \|\bPi_r^\ast \wh\bA_{r_{\max}}\be_j\| 
= \big\| (\bPi_r^\ast - \wh\bA_{r_{\max},1:r} \wh\bA_{r_{\max},1:r}^\top) \wh\bA_{r_{\max}} \be_j \big\|
= \cO_P(r_m) ,
\end{equation}
uniformly over $r < j \leq r_{\max}$,
where the second equality used $\wh\bA_{r_{\max},1:r}^\top \wh\bA_{r_{\max}} \be_j = \zero$,
and the last used $\|\wh\bA_{r_{\max}} \be_j\| = 1$ and \eqref{eqn: dc_top_r_rate}.

Note that by construction, $\bSigma_x \bPi_r^\ast = \bPi_r^\ast \bSigma_x$ and hence $\bs_j^\top \bSigma_x \bn_j =0$.
Also, by construction, the largest eigenvalue of $(\bI_p - \bPi_r^\ast) \bSigma_x (\bI_p - \bPi_r^\ast)$ is exactly $\lambda_{r+1}(\bSigma_x)$.
Hence we have the Rayleigh quotient $\bn_j^\top \bSigma_x \bn_j / \|\bn_j\|^2 \leq \lambda_{r+1}(\bSigma_x) \leq \|\bSigma_\epsilon\| \leq C_\epsilon$ by \eqref{eqn: cov}, Assumptions~\ref{ass: loading} and \ref{ass: noise}
(note that when $\bn_j = \zero$, we have trivially $\bn_j^\top \bSigma_x \bn_j \leq C_\epsilon \cdot \|\bn_j\|^2$, which is sufficient for the below).
Combining these with \eqref{eqn: s_j_rate}, we finally have
\begin{align*}
    \|\bb_j\|^2 &= 
    \big\| \bM^\top \bs_j + \bM^\top \bn_j \big\|^2
    =
    \bs_j^\top \bSigma_x \bs_j 
    + 2 \bs_j^\top \bSigma_x \bn_j 
    + \bn_j^\top \bSigma_x \bn_j \\
    &\leq
    \|\bSigma_x\| \cdot \|\bs_j\|^2
    + \bn_j^\top \bSigma_x \bn_j 
    \leq
    \cO_P(p r_m^2)
    + C_\epsilon \cdot \|\bn_j\|^2
    \leq
    \cO_P(p r_m^2) + C_\epsilon ,
\end{align*}
uniformly over $r < j \leq r_{\max}$,
where the last inequality used the fact that 
$\|\bn_j\|^2 = 1 - \|\bs_j\|^2 \leq 1$ from $1= \|\bs_j\|^2 + \|\bn_j\|^2$.
This thus concludes the proof of this lemma.
\end{proof}

\begin{lemma}[Proposition in Section~2 of \citeauthor{Marsaglia2006}, \citeyear{Marsaglia2006}]
\label{lem: ratio_normal}
For any two jointly normal variables $z$ and $w$ with means $\mu_z$, $\mu_w$, positive variances $\sigma_z^2$, $\sigma_w^2$, and correlation $\rho \in (-1,1)$, the ratio $z/w$ is distributed as
\[
\frac{1}{r}\left( \frac{a+x}{b+y} \right) + s ,
\]
where $x,y$ are independent standard normal, and
\begin{align*}
    r := \frac{\sigma_w}{\sigma_z \sqrt{1-\rho^2}} ,
    \quad
    s := \rho \frac{\sigma_z}{\sigma_w},
    \quad
    b := \frac{\mu_w}{\sigma_w} ,
    \quad
    a := \frac{\mu_z/\sigma_z - \rho \mu_w/\sigma_w}{\sqrt{1-\rho^2}} .
\end{align*}
In particular, if $\mu_z = \mu_w=0$, then $z/w$ is distributed as
\[
\frac{\sigma_z \sqrt{1-\rho^2}}{\sigma_w} c^* + \rho \frac{\sigma_z}{\sigma_w} ,
\]
where $c^* \sim \mathrm{Cauchy}(0,1)$.
\end{lemma}

\begin{proof}[Proof of Lemma~\ref{lem: ratio_normal}]
This is a direct quote from proposition in Section~2 of \cite{Marsaglia2006}, up to rewording and explicitly stating the constants according to the discussion therein.
\end{proof}

\begin{lemma}
\label{lem: median_ratio_quadratic}
Let $\bu$ be uniformly distributed on $\S^{h-1}$ with $h \geq 1$,
and let $\ba, \bb \in \R^{h}$ be deterministic nonzero vectors.
Then it holds that
\begin{equation*}
\mathrm{median}_{\bu \in \S^{h-1}}\left( \frac{ (\ba^{\top} \bu)^2 }{ (\bb^{\top} \bu)^2 } \right)
= \frac{\|\ba\|^2}{\|\bb\|^2} .
\end{equation*}
Moreover, there exists an absolute constant $c_0 >0$ such that, for every $\delta \in (0,1)$,
\begin{align*}
    \P\left(
    \frac{(\ba^\top \bu)^2}{(\bb^\top \bu)^2}
    \leq (1-\delta) \frac{\|\ba\|^2}{\|\bb\|^2} \right)
    &\leq
    \frac{1}{2} - c_0 \delta , \\
    \P\left(
    \frac{(\ba^\top \bu)^2}{(\bb^\top \bu)^2}
    \leq (1+\delta) \frac{\|\ba\|^2}{\|\bb\|^2} \right)
    &\geq
    \frac{1}{2} + c_0 \delta ,
\end{align*}
and these bounds hold uniformly over $h$, $\ba$, and $\bb$.
\end{lemma}

\begin{proof}[Proof of Lemma~\ref{lem: median_ratio_quadratic}]
As $\bu$ is uniformly distributed on the unit sphere $\S^{h-1}$, it can be represented as $\bu = \bz/\|\bz\|$ where $\bz \sim \cN(\zero, \bI_h)$.
Hence we have
\begin{align*}
    \frac{ (\ba^{\top} \bu)^2 }{ (\bb^{\top} \bu)^2 }
    &=
    \frac{ (\ba^{\top} \bz/\|\bz\|)^2 }{ (\bb^{\top} \bz/\|\bz\|)^2 }
    =
    \frac{ (\ba^{\top} \bz)^2 }{ (\bb^{\top} \bz)^2 }
    =
    \frac{\|\ba\|^2}{\|\bb\|^2} \cdot
    \frac{x^2}{y^2} ,
\end{align*}
where $x := \ba^{\top} \bz /\|\ba\| \sim \cN(0, 1)$ and $y := \bb^{\top} \bz /\|\bb\| \sim \cN(0, 1)$.
Denote the correlation between $x$ and $y$ by $\rho := (\ba^\top \bb) /(\|\ba\| \|\bb\|)$.
Note that $|\rho|=1$ if and only if $\ba= c \cdot \bb$ for some constant $c\neq 0$. The desired result holds trivially if $\ba= c \cdot \bb$, and hence it suffices to consider the case where $|\rho|\neq 1$.

Then by Lemma~\ref{lem: ratio_normal}, let $c^*$ be a random sample from standard Cauchy, and we have
\begin{align*}
    &\;\quad
    \P\Bigg( \frac{ (\ba^{\top} \bu)^2 }{ (\bb^{\top} \bu)^2 } \geq \frac{\|\ba\|^2}{\|\bb\|^2} \Bigg)
    =
    \P\Bigg( \frac{x^2}{y^2} \geq 1 \Bigg)
    = 
    \P\Big( 
    \big( \sqrt{1-\rho^2} \cdot c^* + \rho \big)^2
    \geq 1 \Big) \\
    &=
    \P\Big( 
    \sqrt{1-\rho^2} \cdot c^* + \rho 
    \geq 1 \Big)
    + \P\Big( 
    \sqrt{1-\rho^2} \cdot c^* + \rho 
    \leq -1 \Big) \\
    &=
    \P\Bigg( c^* \geq \sqrt{\frac{1-\rho}{1+\rho}} \Bigg)
    + \P\Bigg( c^* \leq - \sqrt{\frac{1+\rho}{1-\rho}} \Bigg)
    =
    \P\Bigg( c^* \geq \sqrt{\frac{1-\rho}{1+\rho}} \Bigg)
    + \P\Bigg( c^* \geq \sqrt{\frac{1+\rho}{1-\rho}} \Bigg) \\
    &=
    \frac{1}{2} - \frac{1}{\pi} \mathrm{arctan}\Bigg( \sqrt{\frac{1-\rho}{1+\rho}} \Bigg)
    + \frac{1}{2} - \frac{1}{\pi} \mathrm{arctan}\Bigg( \sqrt{\frac{1+\rho}{1-\rho}} \Bigg)
    = 
    1 - \frac{1}{\pi} \cdot \frac{\pi}{2}
    = \frac{1}{2} ,
\end{align*}
as desired, where the last line used the fact that $\mathrm{arctan}(x) + \mathrm{arctan}(1/x) = \pi/2$ for any $x>0$.
This shows the desired median.

For every $s > 0$, by Lemma~\ref{lem: ratio_normal} as above again, 
and by direct computation, we obtain
\begin{align*}
    \P\Bigg( \frac{x^2}{y^2} \leq s \Bigg) &=
    \frac{1}{\pi}
    \Bigg\{
    \arctan\Bigg( \frac{\sqrt{s} - \rho}{\sqrt{1-\rho^2}} \Bigg)
    - \arctan\Bigg( \frac{-\sqrt{s} - \rho}{\sqrt{1-\rho^2}} \Bigg)
    \Bigg\} \\
    &=
    \frac{1}{2}
    +
    \frac{1}{\pi}
    \arctan\Bigg( \frac{s-1}{2 \sqrt{s} \sqrt{1-\rho^2}} \Bigg) ,
\end{align*}
so that for $\delta \in (0,1)$,
using $\sqrt{1-\rho^2} \leq 1$, we have
\begin{align*}
    \P\Bigg( \frac{x^2}{y^2} \leq 1+\delta \Bigg)
    - \frac{1}{2}
    &=
    \frac{1}{\pi}
    \arctan\Bigg( \frac{\delta}{2 \sqrt{1+\delta} \sqrt{1-\rho^2}} \Bigg) \\
    &\geq
    \frac{1}{\pi}
    \arctan\Bigg( \frac{\delta}{2 \sqrt{2}} \Bigg)
    \geq
    \frac{\delta}{4\pi \sqrt{2}} ,
\end{align*}
and similarly,
\begin{align*}
    \frac{1}{2}
    - \P\Bigg( \frac{x^2}{y^2} \leq 1-\delta \Bigg)
    &=
    \frac{1}{\pi}
    \arctan\Bigg( \frac{\delta}{2 \sqrt{1-\delta} \sqrt{1-\rho^2}} \Bigg) \\
    &\geq
    \frac{1}{\pi}
    \arctan\Bigg( \frac{\delta}{2 \sqrt{2}} \Bigg)
    \geq
    \frac{\delta}{4\pi \sqrt{2}} ,
\end{align*}
where the last line used $\arctan(v) \geq v/2$ for $v \in [0,1]$.
Thus we may take $c_0 = 1/(4\pi \sqrt{2})$, which completes the proof of this lemma.
\end{proof}

\subsection{Proofs of corollaries and propositions}
\label{subsec: proof_cor}

\begin{proof}[Proof of Corollary~\ref{cor: dc_star}]
Let $\wh\bA_r$ denote the DC estimator computed from all observations with the deterministic block size $b = r$.
Note that the resulting number of blocks satisfies $\floor{T/r} \asymp T$ as $r$ is fixed,
$\overline{\cT} = o(T)$ by \eqref{eqn: thm_r_contam_condition},
and hence by Theorem~\ref{thm: DC_ano_clean} we have
\begin{equation}
\label{eqn: dc_rate_T}
d(\wh\bA_r, \bA) =
\cO_P\left(
\frac{1}{\sqrt{T}}
+ \frac{1}{p}
+ \frac{\overline{\cT}}{T}
\right) .
\end{equation}
On the event $\{\wh{r} = r\}$, the estimator $\wh\bA_{\wh{r}} = \wh\bA_r$. Hence for every deterministic $u > 0$,
\[
\P\left( d(\wh\bA_{\wh{r}}, \bA) > u \right)
\leq
\P(\wh{r} \neq r)
+
\P\left( d(\wh\bA_r, \bA) > u \right) ,
\]
which, together with Theorem~\ref{thm: r_hat_consistency} and \eqref{eqn: dc_rate_T}, concludes this proof.
\end{proof}

\begin{proof}[Proof of Corollary~\ref{cor: dc_fac_common}]
Fix any $t\in[T]$,
recall that $\cT=\emptyset$ here,
and recall $\bv_t = \bu_t / \|\bu_t\|$.
Since $\bv_t$ is uniform on $\S^{r+p-1}$, we have $\E[\bv_t \bv_t^\top] = (r+p)^{-1} \bI_{r+p}$.
Then by \eqref{eqn: elliptical_obs}, 
\eqref{eqn: xt_ellip_form},
and $\zeta_t > 0$ almost surely, 
we obtain
\begin{align*}
    \E\left[ \frac{r+p}{p \zeta_t^2} \|\bx_t\|^2 \right]
    &=
    \frac{r+p}{p}
    \E[\|\bM \bv_t\|^2] 
    =
    \frac{1}{p} \tr(\bM \bM^\top) 
    =
    \frac{1}{p} \tr(\bSigma_x) 
    \leq
    r C_\lambda + C_\epsilon , 
    \\
    \E\left[ \frac{r+p}{p \zeta_t^2} \|\bepsilon_t\|^2 \right]
    &=
    \frac{1}{p} \tr(\bSigma_\epsilon)
    \leq
    C_\epsilon , 
    \\
    \E\left[ \frac{r+p}{\zeta_t^2}
    \|\ddot\bA \ddot\bA^\top \bepsilon_t\|^2 \right]
    &=
    \tr\big(
    \ddot\bA \ddot\bA^\top
    \bSigma_\epsilon
    \big)
    \leq
    r C_\epsilon ,
\end{align*}
where, particularly, note that these expectations are finite without any moment condition on $\zeta_t$, since each $\zeta_t^2$ is canceled out by itself within $\|\bx_t\|^2$, $\|\bepsilon_t\|^2$, and $\|\ddot\bA \ddot\bA^\top \bepsilon_t\|^2$, respectively.
Hence by $\zeta_t^2 /p = \cO_P(1)$, it holds that
\begin{equation}
\label{eqn: fac_common_step_moment}
\frac{\|\bx_t\|^2}{p} = \cO_P(1) ,
\qquad
\frac{\|\bepsilon_t\|^2}{p} = \cO_P(1) ,
\qquad
\|\ddot\bA \ddot\bA^\top \bepsilon_t\|^2 = \cO_P(1) .
\end{equation}
Hence by model~\eqref{eqn: vfm} and definitions of the estimators and $\wh\bH$, we obtain
\begin{align*}
    &\;\quad
    \big\| \wh\bf_t - \wh\bH \bf_t \big\|^2 = 
    \left\| \frac{1}{\sqrt{p}}
    \wh\bA_{\wh{r}}^\top \bepsilon_t \right\|^2
    =
    \frac{1}{p} \big\| \wh\bA_{\wh{r}} \wh\bA_{\wh{r}}^\top \bepsilon_t \big\|^2 \\
    &\leq
    \frac{2}{p} \big\| \ddot\bA \ddot\bA^\top \bepsilon_t \big\|^2
    + \frac{2}{p} 
    \big\| \wh\bA_{\wh{r}} \wh\bA_{\wh{r}}^\top - \ddot\bA \ddot\bA^\top \big\|^2 
    \cdot \big\| \bepsilon_t \big\|^2 \\
    &=
    \cO_P\left( 
    \frac{1}{p}
    + \big\| \wh\bA_{\wh{r}} \wh\bA_{\wh{r}}^\top - \ddot\bA \ddot\bA^\top \big\|^2 
    \right)
    =
    \cO_P\left(
    \frac{1}{T}
    +\frac{1}{p} 
    \right) ,
\end{align*}
as desired, where the last line used \eqref{eqn: fac_common_step_moment} and Corollary~\ref{cor: dc_star}.

Similarly for the common component,
as $\wh\bchi_t - \bchi_t = \big( \wh\bA_{\wh{r}} \wh\bA_{\wh{r}}^\top - \ddot\bA \ddot\bA^\top \big) \bx_t + \ddot\bA \ddot\bA^\top \bepsilon_t$,
by \eqref{eqn: fac_common_step_moment} and Corollary~\ref{cor: dc_star} again, we have
\begin{align*}
    \frac{1}{p} \big\| \wh\bchi_t - \bchi_t \big\|^2 &\leq
    2 \big\| \wh\bA_{\wh{r}} \wh\bA_{\wh{r}}^\top - \ddot\bA \ddot\bA^\top \big\|^2 \frac{\|\bx_t\|^2}{p}
    + \frac{2}{p} \big\| \ddot\bA \ddot\bA^\top \bepsilon_t \big\|^2
    =
    \cO_P\left(
    \frac{1}{T}
    +\frac{1}{p} 
    \right) ,
\end{align*}
as desired.
It remains to show $\wh\bH$ is well conditioned.
As $\P(\wh{r} = r) \to 1$ by Theorem~\ref{thm: r_hat_consistency}, it suffices to work under the event $\{\wh{r} = r\}$.
Then $\wh\bH$ being well conditioned is straightforward by 
\eqref{ass: loading_star}
and noting that
\begin{align*}
    &\;\quad
    \left\| \frac{1}{p} \bA^\top \bA - \wh\bH^\top \wh\bH \right\|
    =
    \left\| \frac{1}{p} \bA^\top \bA - \frac{1}{p} \bA^\top \wh\bA_{\wh{r}} \wh\bA_{\wh{r}}^\top \bA \right\|
    =
    \frac{1}{p} \big\| \bA^\top (\bI_p - \wh\bA_{\wh{r}} \wh\bA_{\wh{r}}^\top) \bA \big\| \\
    &=
    \frac{1}{p} \big\| (\bI_p - \wh\bA_{\wh{r}} \wh\bA_{\wh{r}}^\top) \bA \big\|^2
    \leq
    \frac{1}{p} \|\bA\|^2 
    \cdot \big\| \wh\bA_{\wh{r}} \wh\bA_{\wh{r}}^\top - \ddot\bA \ddot\bA^\top \big\|^2 
    = o_P(1) ,
\end{align*}
where the third equality used the fact that $(\bI_p - \wh\bA_{\wh{r}} \wh\bA_{\wh{r}}^\top)$ is an orthogonal projection matrix, 
the inequality used $\ddot\bA \ddot\bA^\top \bA = \bA$ by construction,
and the last equality used Corollary~\ref{cor: dc_star} again.
This thus concludes this proof.
\end{proof}

\begin{proof}[Proof of Proposition~\ref{prop: singleton_failure}]
As in the proof of Theorem~\ref{thm: DC_ano_clean},
recall the latent observation $\bx_t^0$.
It is direct to check Assumption~\ref{ass: elliptical} is satisfied.
By construction in the statement, the eigenvalues of $\bA^\top \bA$ are $\kappa_p(p + p^\delta)$ and $p^\delta$, respectively,
so Assumption~\ref{ass: loading} holds with $c_\lambda = 1$.
Assumption~\ref{ass: noise} holds with $c_\epsilon = C_\epsilon = 1$.
Also, $|\cT| = \overline{\cT} = 1$ and satisfies Assumption~\ref{ass: contamination}.

Next, note that we can represent
\[
\bx_t^0 = a_p G_t \be_1 + \by_t ,
\qquad
a_p^2 := \kappa_p (p+ p^\delta) +1 ,
\qquad
\be_1^\top \by_t =0 ,
\]
where $G_t \sim \cN(0,1)$ independent of $\by_t$, 
and $\E[\|\by_t\|^2] = p^\delta + p-1 \leq p + p^\delta$.
Hence with $\bQ_t^0 := \bx_t^0 (\bx_t^0)^\top /\|\bx_t^0\|^2$, we have
\begin{equation}
\label{eqn: Qt0_ee}
\big\| \bQ_t^0 - \be_1\be_1^\top \big\|_F
=
\sqrt{2} \cdot \frac{\|\by_t\|}{\sqrt{a_p^2 G_t^2 + \|\by_t\|^2}}
\leq
\sqrt{2} \cdot \min\left\{ 1, \frac{\|\by_t\|}{a_p \cdot |G_t|} \right\} .
\end{equation}

Note that conditional on $G_t$, by the Jensen's inequality, we have
\[
\E\left[
\min\left\{ 1, \frac{\|\by_t\|}{a_p \cdot |G_t|} \right\}
\mid G_t \right]
\leq
\min\left\{ 1, \frac{\sqrt{ \E[\|\by_t\|^2] }}{a_p \cdot |G_t|} 
\right\}
\leq
\min\left\{ 1, \frac{\kappa_p^{-1/2}}{|G_t|}
\right\} ,
\]
which, together with the notation $\phi(\cdot)$ representing the density function of $\cN(0,1)$, 
using $\kappa_p^{-1/2} < 1$, implies
\begin{align*}
    \E\left[ 
    \min\left\{ 1, \frac{\kappa_p^{-1/2}}{|G_t|} \right\}
    \right]
    &=
    \P\big( |G_t| \leq \kappa_p^{-1/2} \big) 
    + 2 \kappa_p^{-1/2} \int_{\kappa_p^{-1/2}}^\infty \frac{\phi(v)}{v} \, \mathrm{d}v \\
    &\leq
    2 \phi(0) \kappa_p^{-1/2} 
    + 2 \phi(0) \kappa_p^{-1/2} \int_{\kappa_p^{-1/2}}^1 \frac{1}{v} \, \mathrm{d}v
    + 
    2 \kappa_p^{-1/2} \int_1^\infty \phi(v) \, \mathrm{d}v \\
    &\leq
    C \kappa_p^{-1/2} \{1 + \log(1/ \kappa_p^{-1/2}) \}
    \leq
    C \frac{\log(\kappa_p)}{\sqrt{\kappa_p}} ,
\end{align*}
which, combined with \eqref{eqn: Qt0_ee}, gives
\begin{equation}
\label{eqn: singleton_clean_bound}
\E\big[ \|\bQ_t^0 -\be_1 \be_1^\top \|_F \big]
\leq
C \frac{\log(\kappa_p)}{\sqrt{\kappa_p}} .
\end{equation}

Let $J$ be the anomaly location,
and note that $\bw := \bz/\|\bz\|$ is independent and uniform on the unit sphere in $\R^p$.
Then define the matrix
\[
\bL := 
\frac{T-1}{T} \be_1 \be_1^\top 
+ \frac{1}{T} \bw\bw^\top .
\]
By the triangle inequality, \eqref{eqn: singleton_clean_bound},
and the independence of $J$ from the clean data, we have
\begin{align*}
    \E\big[ \|\wh\bP^{(1)} - \bL\| \big]
    &\leq
    \E\Bigg[\frac{1}{T} \sum_{t \in [T] \setminus \{J\} } \big\| \bQ_t^0 -\be_1 \be_1^\top \big\|_F \Bigg] \\
    &=
    \frac{T-1}{T} \cdot
    \E\big[ \big\| \bQ_1^0 -\be_1\be_1^\top \big\|_F \big]
    \leq
    C \frac{\log(\kappa_p)}{\sqrt{\kappa_p}} ,
\end{align*}
which, under \eqref{eqn: singleton_growth}, implies $\|\wh\bP^{(1)} - \bL\| = o_P(1/T)$.

Since $\E[(\be_1^\top\bw)^2] = 1/p$, we have $\be_1^\top\bw = o_P(1)$.
Note that the range of $\bL$ is $\spa(\be_1, \bw)$, and its two positive eigenvalues have sum $1$ and product $\{(T - 1)/T^2\} (1 - (\be_1^\top\bw)^2)$.
As $\lambda_1(\bL) \leq 1$, we hence have
\[
\lambda_2(\bL)
=
\frac{(T-1) (1 - (\be_1^\top\bw)^2)}{T^2 \lambda_1(\bL)}
\geq
\frac{(T-1) (1 - (\be_1^\top\bw)^2)}{T^2}
\geq
\frac{1}{4T} ,
\]
with probability approaching one, while $\lambda_3(\bL) = 0$.
Hence with the definition
\[
\bv := \frac{\bw - \be_1^\top\bw \be_1}{\sqrt{1- (\be_1^\top\bw)^2}},
\qquad
\bQ := \be_1\be_1^\top + \bv\bv^\top ,
\]
by the Davis--Kahan theorem \citep[Theorem~2]{yu2015useful}, 
we have
\[
\big\| \wh\bA^{(1)}(\wh\bA^{(1)})^\top - \bQ \big\|_F
\leq 
C T \big\|\wh\bP^{(1)}- \bL \big\|
= o_P(1) .
\]
Finally, from $\E[(\be_2^\top\bw)^2] = 1/p$ and $(\be_1^\top\bw)^2 =o_P(1)$, we have
\[
\big\| \bQ - \big(\be_1\be_1^\top + \be_2\be_2^\top \big) \big\|_F^2
=
2 \big\{ 1 - (\be_2^\top\bv)^2 \big\}
, \qquad
(\be_2^\top\bv)^2 =
\frac{(\be_2^\top \bw)^2}{ 1- (\be_1^\top\bw)^2}
= o_P(1) ,
\]
so the reverse triangle inequality now gives $d(\wh\bA^{(1)},\bA) - \sqrt{2} = o_P(1)$ as desired, which concludes this proof.
\end{proof}

\clearpage
\section{Additional discussion of Remark~\ref{remark: spatial_sign_cov}}
\label{app:rem:spatial}

We discuss here the performance of the DC estimator applied with the block size $b=1$ in~\eqref{eqn: block_cov_vfm}, which fails under our general setup.
In particular, we present the following proposition, with its proof in Section~\ref{subsec: proof_cor}.

\begin{proposition}
\label{prop: singleton_failure}
Fix any $\delta > 0$, let $r = 2$,
and let $\be_1, \be_2$ be the first two standard basis vectors in $\R^p$. Then with $\kappa_p \geq 2$, construct the factor loading matrix
\begin{equation*}
\bA =
\left[
\sqrt{\kappa_p (p + p^\delta)} \, \be_1 ,
p^{\delta/2}\be_2
\right] .
\end{equation*}
Let $\bff_t$ and $\bepsilon_t$ be independent standard Gaussian vectors that are independently across observations.
Also, replace one uniformly selected observation $\bx_t$ by an anomaly vector $\bz \sim \cN(\zero, \bI_p)$,
where $\bz$ is independent of all latent clean observations,
and its location is independent of all latent clean observations and of $\bz$.
Suppose
\begin{equation}
\label{eqn: singleton_growth}
\frac{1}{p} , \;\;
\frac{p^{2 \max\{1 - \delta,0\}}}{T} , \;\;
\frac{T \log(\kappa_p)}{\sqrt{\kappa_p}} 
= o(1) .
\end{equation}
Then the above setup satisfies Assumptions~\ref{ass: elliptical}, \ref{ass: loading}, \ref{ass: noise}, and \ref{ass: contamination}, 
with $\overline{\cT} = 1$.
Then with the notations
\begin{equation}
\label{eqn: singleton_estimator}
\wh\bP^{(1)} := 
\frac{1}{T} \sum_{t=1}^T
\frac{\bx_t\bx_t^\top}{\|\bx_t\|^2},
\qquad
\wh\bA^{(1)} := 
\mathrm{SVD}_2(\wh\bP^{(1)}) ,
\end{equation}
we have $d(\wh\bA^{(1)},\bA) - \sqrt{2} = o_P(1)$.
\end{proposition}

The conditions in~\eqref{eqn: singleton_growth} are compatible for every fixed $\delta > 0$,
such as $T \asymp p^a$ for any fixed $a > 2 \max\{1 - \delta,0\}$, and $\kappa_p \asymp T^3$.
Proposition~\ref{prop: singleton_failure} establishes a failure of using block size $b=1$ in~\eqref{eqn: block_cov_vfm}.
In contrast, by our construction with $b= r= 2$, Theorem~\ref{thm: DC_ano_clean} gives
\[
d(\wh\bA_2,\bA)
= 
\cO_P\left(
\frac{p^{\max\{1 - \delta,0\}}}{\sqrt{T}}
+ p^{-\delta}
\right)
= o_P(1) .
\]

\clearpage
\section{Additional numerical results}
\label{add:sim}

\subsection{Main numerical experiments under various dimensions}

We report here the results as in the experiments corresponding to Figure~\ref{fig: sim_nu_epsilon}, but under dimensions $(T,p) \in \{(50,50), (200,200)\}$.
The results are shown in Figures~\ref{fig: sim_nu_epsilon_(50,50)}
and \ref{fig: sim_nu_epsilon_(200,200)}, respectively.
From these figures, the discussions to Figure~\ref{fig: sim_nu_epsilon} remain valid.

\begin{figure}[ht!]
\centering
\includegraphics[width=.85\linewidth]{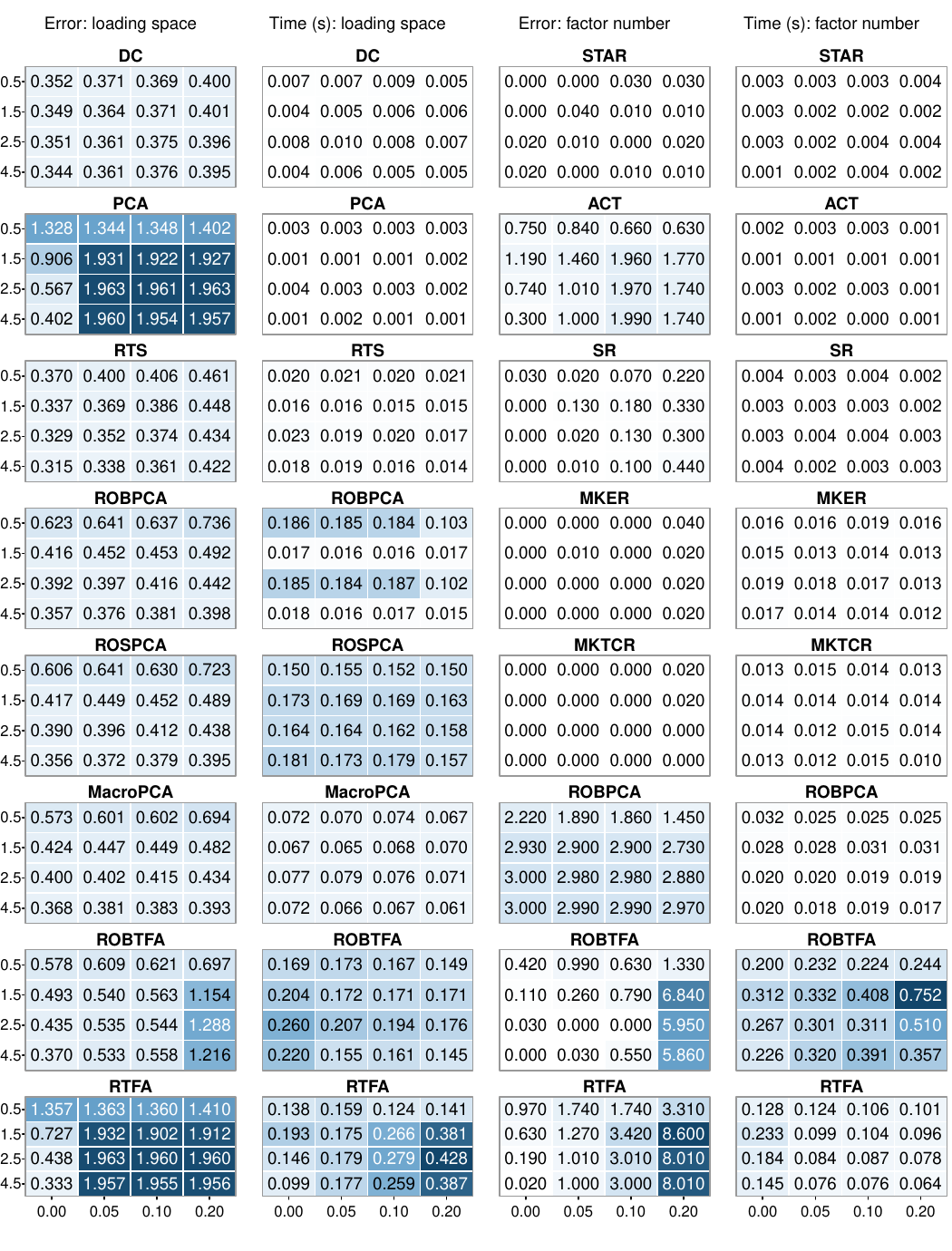}
{\setstretch{1}
\caption{
Performance of loading space and factor number estimation for $(T,p)=(50,50)$.
Within each heatmap, rows correspond to
$\nu \in \{ 0.5, 1.5, 2.5, 4.5 \}$
and columns correspond to
$\varepsilon \in \{ 0, 0.05, 0.1, 0.2 \}$.
Darker shading indicates larger values within each column.
}
\label{fig: sim_nu_epsilon_(50,50)}}
\end{figure}

\begin{figure}[ht!]
\centering
\includegraphics[width=.85\linewidth]{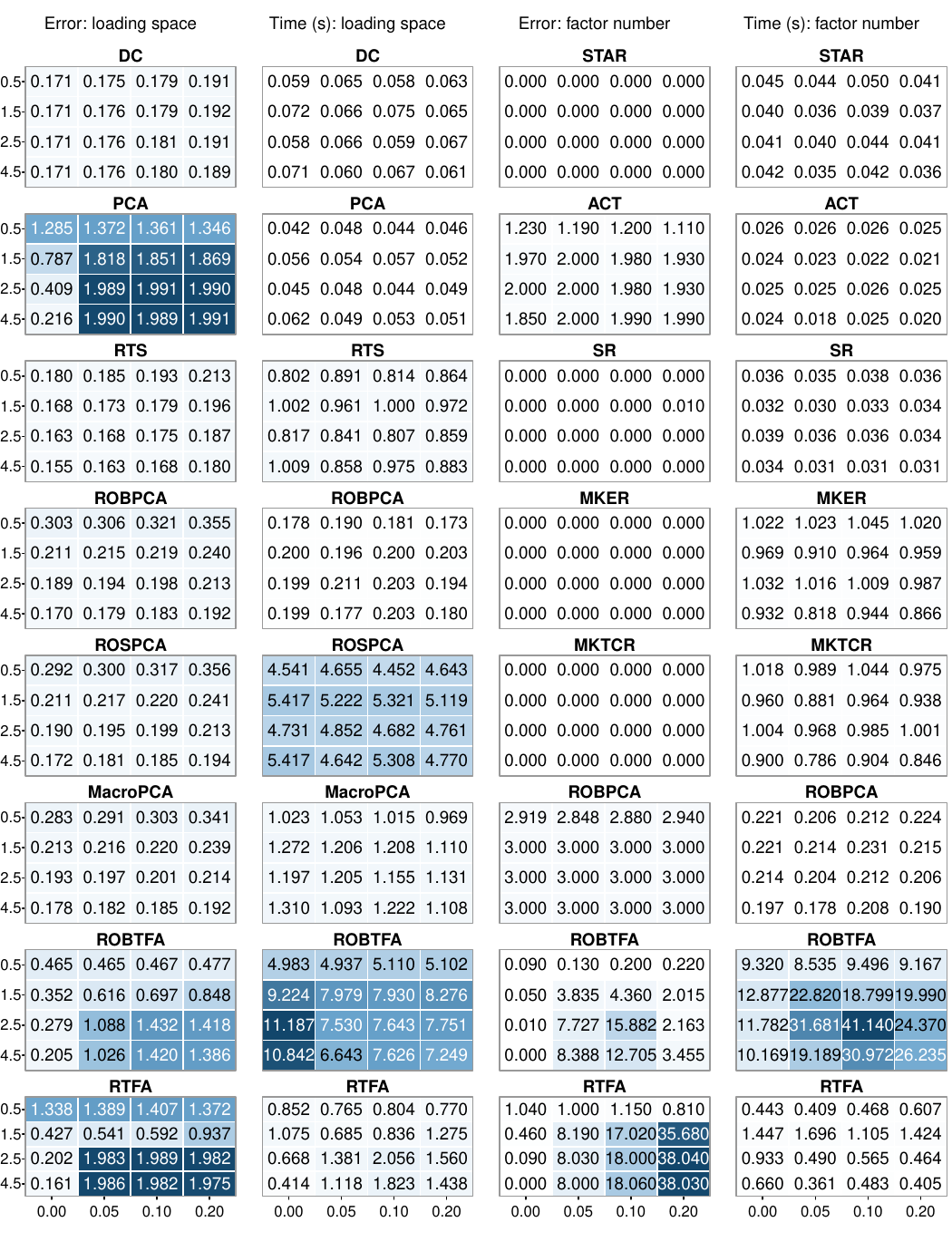}
{\setstretch{1}
\caption{
Performance of loading space and factor number estimation for $(T,p)=(200,200)$.
Within each heatmap, rows correspond to
$\nu \in \{ 0.5, 1.5, 2.5, 4.5 \}$
and columns correspond to
$\varepsilon \in \{ 0, 0.05, 0.1, 0.2 \}$.
Darker shading indicates larger values within each column.
}
\label{fig: sim_nu_epsilon_(200,200)}}
\end{figure}

\subsection{Common component estimation}
\label{subsec: sim_common_component}

With Corollary~\ref{cor: dc_fac_common},
we further investigate common component estimation under heavy tails and anomalies, using the data-generating process in Section~\ref{subsec: sim}.
We fix $r = 2$ and $a = \phi = 0$, and consider $(T,p) \in \{(50,50), (100,100), (200,200), (400,400)\}$.
For each dimension pair, we vary the degrees of freedom over $\nu \in \{0.5, 1.5, 2.5, 4.5\}$ and the fraction of anomalies over $\varepsilon \in \{0, 0.05, 0.1, 0.2\}$.
The remaining specifications follow the main experiment in Section~\ref{subsec: sim}.
Each configuration is repeated $100$ times, with all methods applied to the same observations within each replication.

We compare DC with the seven loading space estimators
considered in Section~\ref{subsec: sim}, namely
\texttt{PCA}, 
\texttt{RTS}, 
\texttt{ROBPCA},
\texttt{ROSPCA}, 
\texttt{MacroPCA}, 
\texttt{ROBTFA},
and \texttt{RTFA}.
The true factor number $r$ is provided to every method, and hence DC uses $b = r = 2$.
For each estimated loading space, we form its orthogonal projector and estimate the common component by projecting the observed $\bx_t$ onto that space, as in~\eqref{eqn: dc_fac_common_est}.
All observations, including anomalies, are used to estimate the loading space.
For evaluation, we retain the latent common components generated before any observations are replaced and measure reconstruction error at the clean time points.
Specifically, for each replication we compute
\begin{equation}
\label{eqn: sim_common_relative_error}
\left(
\frac{\sum_{t\in [T] \setminus \cT} \big\| \wh\bchi_t - \bchi_t \big\|^2}{\sum_{t\in [T] \setminus \cT}
\big\| \bchi_t \big\|^2}
\right)^{1/2} .
\end{equation}
It is noteworthy that the anomaly indices are only used in evaluation, rather than provided to any estimator.
The relative normalization in~\eqref{eqn: sim_common_relative_error} facilitates comparisons across dimensions and tail parameters, and similar metrics are used in e.g.\ \cite{CenLam2025}.

\begin{figure}[ht!]
\centering
\includegraphics[width=.85\linewidth]{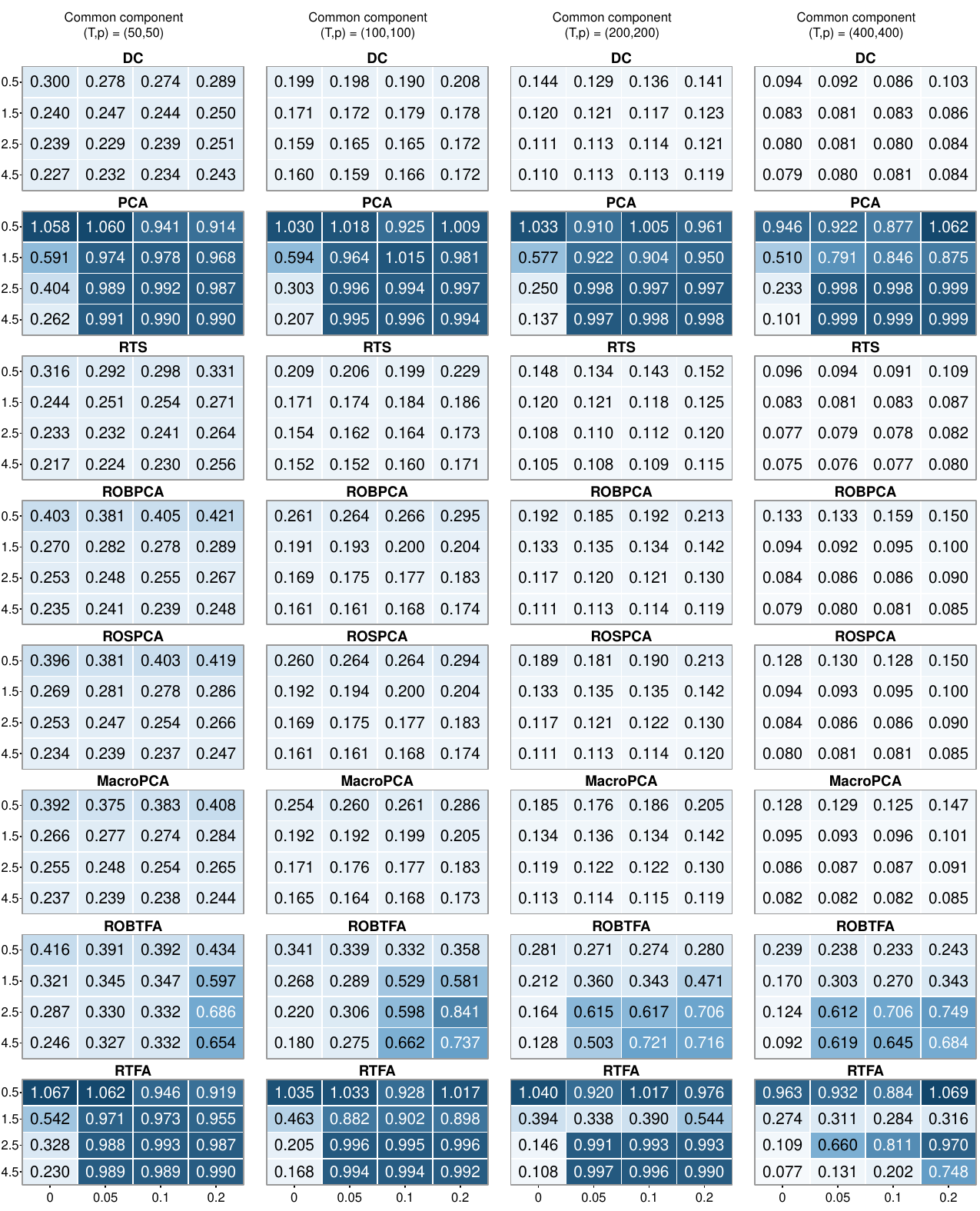}
{\setstretch{1}
\caption{
Common component estimation under heavy tails and anomalies.
Within each heatmap, rows correspond to
$\nu \in \{0.5, 1.5, 2.5, 4.5\}$ and columns correspond to
$\varepsilon \in \{0, 0.05, 0.1, 0.2\}$.
Cell values report the average relative common component error in~\eqref{eqn: sim_common_relative_error} over $100$ replications.
Darker shading indicates larger values within each column.
}
\label{fig: sim_common_component}}
\end{figure}

\clearpage 
\section{Additional information for the real data analysis}
\label{subsec: real_data_details}

The raw observation at each timestamp consists of $10 \times 10$ categories with each forming one portfolio, where the stocks are weighted by their market equity.
If no stocks are allocated to a category at a given timestamp, the corresponding return is unavailable.
The missing pattern is irregular across categories, and dropping four categories to obtain a complete data set, or restricting attention to the period after July 2009, would discard a substantial portion of the information. We therefore adopt the imputed panels of \cite{CenLam2025}, with the market factor also removed from the dataset following their procedure.
Each monthly observation as a $10 \times 10$ array is then vectorized, and we analyze the resulting vector $\bx_t \in \R^{100}$.

\begin{figure}[ht!]
\centering
\includegraphics[width = \linewidth]{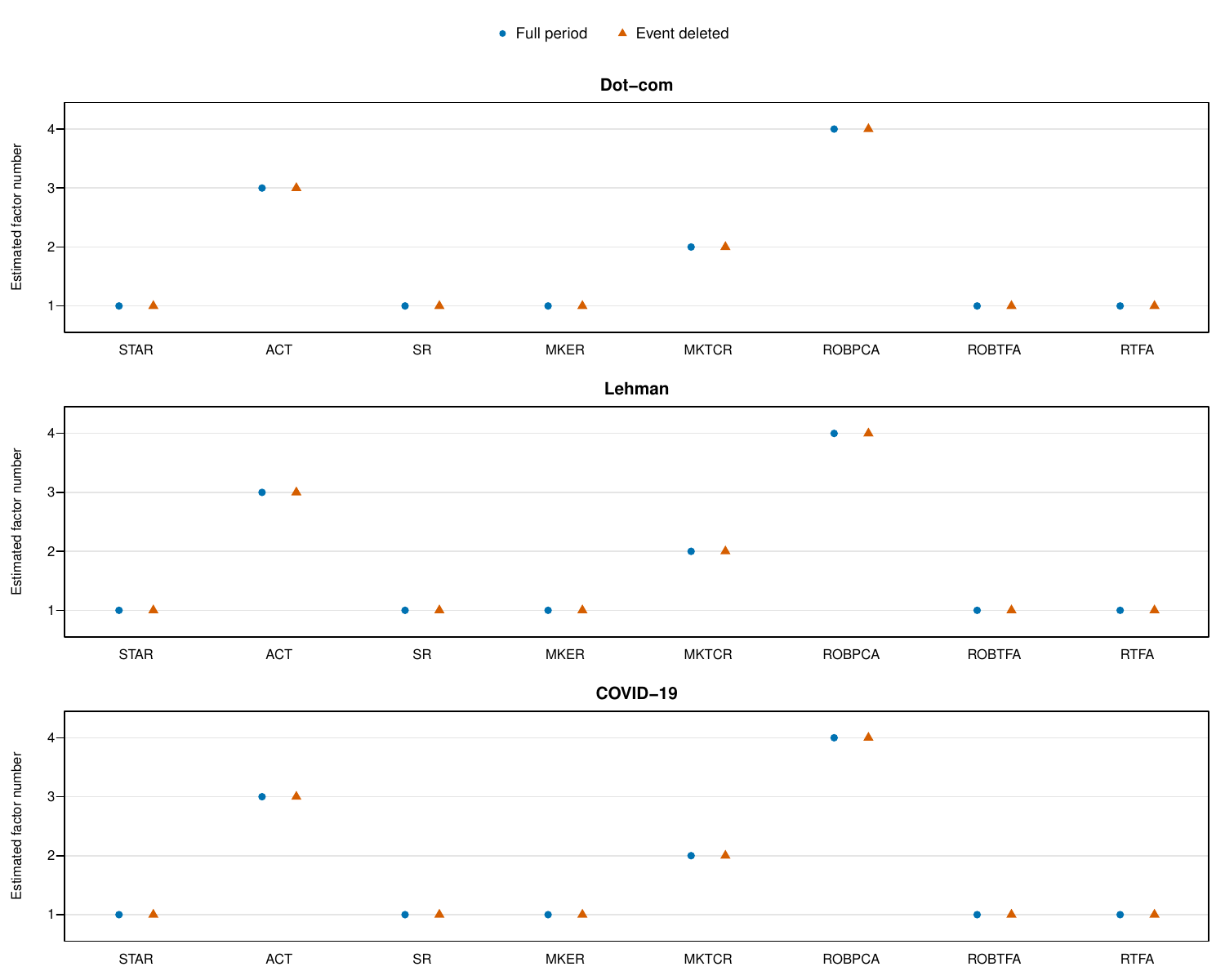}
\caption{
Estimated factor numbers for the complete sample and the samples obtained by separately removing the windows corresponding to the three events.
Blue circles represent complete-sample estimates and orange triangles represent estimates after deletion.
}
\label{fig: real_rank_stability}
\end{figure}

Figure~\ref{fig: real_rank_stability} reports the estimated factor numbers.
Every method gives the same estimate before and after each exclusion, although the selected factor numbers differ across methods.

\end{document}